\documentclass[12pt]{amsart}

\input xy
\xyoption{all}

\usepackage{geometry}
\usepackage{amsmath}
\usepackage{amssymb}
\usepackage{cite}
\usepackage{amsthm}  % see geometry.pdf on how to lay out the page. There's lots.
\newtheorem{ntheorem}{Theorem}[section]
\newtheorem{ncor}[ntheorem]{Corollary}
\newtheorem{nlemma}[ntheorem]{Lemma}
\newtheorem{nprop}[ntheorem]{Proposition}
\newtheorem{ndefinition}[ntheorem]{Definition}
\newtheorem{nclaim}[ntheorem]{Claim}
\newtheorem{nfact}[ntheorem]{Fact}
\newtheorem{hyp}[ntheorem]{Hypothesis}

\newbox\noforkbox \newdimen\forklinewidth
\setbox0\hbox{$\textstyle\smile$}
\setbox1\hbox to \wd0{\hfil\vrule width \forklinewidth depth-2pt
 height 10pt \hfil}
\setbox\noforkbox\hbox{\lower 2pt\box1\lower 2pt\box0\relax}
\def\unionstick{\mathop{\copy\noforkbox}\limits}

\def\nonfork_#1{\unionstick_{\textstyle #1}}

\setbox0\hbox{$\textstyle\smile$}
\setbox1\hbox to \wd0{\hfil{\sl /\/}\hfil}
\setbox2\hbox to \wd0{\hfil\vrule height 10pt depth -2pt width
              \forklinewidth\hfil}
\newbox\doesforkbox
\setbox\doesforkbox\hbox{\lower 2pt\box1 \lower 2pt\box2\lower2pt\box0\relax}
\def\nunionstick{\mathop{\copy\doesforkbox}\limits}

\def\fork_#1{\nunionstick_{\textstyle #1}}

\newcommand{\dnf}{\unionstick}
\newcommand{\dnfty}{\dnf_{\geq \s}}

\newcommand{\cf}{\text{cf }}
\newcommand{\rest}{\upharpoonright}
\newcommand{\id}{\textrm{id}}

\newcommand{\realize}{\vDash}

\newcommand{\s}{\mathfrak{s}}
\newcommand{\bs}{bs}%\te{bs}}

\newcommand{\seq}[1]{\langle #1 \rangle}
\newcommand{\te}[1]{\textrm{#1}}

\newcommand{\gS}{ \textbf{S} }

\title{Tameness and Extending Frames}
\date{May 8, 2014\\
AMS 2010 Subject Classification: 03C48, 03C45.} % delete this line to display the current date

\author{Will Boney}
\email{wboney@cmu.edu}
\address{Department of Mathematical Sciences \\ Carnegie Mellon University \\ Pittsburgh, Pennsylvania, USA}

\begin{document}

\maketitle

%%%%%%%%%%%%%%%%%%%%%%%%%%%%%%%%%%%%%%%%%%%%%%%
\begin{abstract}
We combine two notions in AECs, tameness and good $\lambda$-frames, and show that they together give a very well-behaved nonforking notion in all cardinalities.  This helps to fill a longstanding gap in classification theory of tame AECs and increases the applicability of frames.  Along the way, we prove a complete stability transfer theorem and uniqueness of limit models in these AECs.
\end{abstract}%%%%%%%%%%%%%%%%%%%%%%%%%%%%%%%%%%%%%%%%%%%%%%%

\tableofcontents

%%%%%%%%%%%%%%%%%%%%%%%%%%%%%%%%%%%%%%%
\section{Introduction} \label{introsection}
%%%%%%%%%%%%%%%%%%%%%%%%%%%%%%%%%%%%%%%

In this paper, we combine two recent developments in Abstract Elementary Classes (AECs): tameness and good $\lambda$-frames.  Tameness is a locality condition for Galois types and good $\lambda$-frames are an axiomatic independence notion for $K_\lambda$.  Doing so allows us to extend the good $\lambda$-frame $\s$, which operates only on $\lambda$-sized models, to a good frame $\geq \s$ that is a forking notion for the entire class.  Precisely, we prove the following.

\begin{ntheorem}
If $K$ is $\lambda$-tame for 1- and 2- types, $\s$ is a good $\lambda$-frame, and $K$ satisfies the amalgamation property, then $\geq \s$ is a good frame.  In particular, $K_{\geq \lambda}$ has no maximal models, is stable in all cardinals, and has a unique limit model in each cardinal.
\end{ntheorem}

We define these notions in the next sections, but give some background here.

Tameness is a locality property for Galois types in AECs.  Recall that the syntactic definition of type is not useful in a general AEC, so Galois types were introduced in Shelah \cite{sh300} as a replacement.  Because we will only use Galois types in this paper, we use `type' to mean Galois type.  Tameness was first isolated in Grossberg and VanDieren \cite{tamenessone}, which came from the latter's thesis, and says that two different types over a large model must differ over some smaller model.  Tameness has been used successfully in categoricity transfers (see Grossberg and VanDieren \cite{tamenesstwo} and \cite{tamenessthree} and Lessman \cite{lessmancat}) and stability transfer (see Grossberg and VanDieren \cite{tamenessone}; Baldwin, Kueker, and VanDieren \cite{bkv}; and Lieberman \cite{liebermanrank}).  Unfortunately, not all AECs are tame as Baldwin and Shelah \cite{nonlocality} have constructed an AEC that is not tame from the exact sequences of an almost free, non-Whitehead group which exists in ZFC at $\aleph_1$ and consistently exists in all cardinals.  On the other hand, the author has shown in \cite{tamelc} that tameness follows for all AECs from large cardinals.
\begin{ntheorem}[\cite{tamelc}]
\begin{itemize}
\item If $K$ is an AEC with $LS(K) < \kappa$ and $\kappa$ is strongly compact, then $K$ is $<\kappa$ tame.
\item Suppose there is a proper class of strongly compact cardinals.  Then every AEC is tame.
\end{itemize}
\end{ntheorem}
Additionally, the author and Grossberg have shown in \cite{shorttamedep} that tameness follows from a strong enough independence relation, a phenomenon first observed in \cite{tamenessthree}.

Frames are a notion of nonforking for a general AEC.  They were first defined axiomatically in Shelah \cite{sh600}, which is published as \cite{shelahaecbook}.II.  \cite{sh600} draws on the results of Shelah \cite{sh576} to provide a general example of a frame from categoricity in two consecutive cardinals, a medium number of models in the third, and some non-ZFC axioms; see Theorem \ref{genframeext} for the precise statement.  The first volume of \cite{shelahaecbook}, Jarden and Shelah \cite{jrsh875}, and Jarden and Sitton \cite{jardensitton} are focused on using frames to develop classification theory for AECs.  This is done by taking good $\lambda$-frames and shrinking the class as the size of the models goes up.  We avoid this very complicated process by the use of tameness.  Shelah defines a more general notion of an extended frame $\geq \s$, but does so only as ``an exercise to familiarize the reader with $\lambda$ frames'' \cite{shelahaecbook}(p. 264).  He shows that some of the frame properties follow (see Theorem \ref{basicextension}).  Here we use tameness to derive the remaining properties.  
 Note that we use the definition of frames from the more recent \cite{jrsh875}.  This definition leaves out some of the redundant clauses and, more significantly, does not require the existence of a superlimit model.

Prior to this paper, there has been no work examining frames and tameness together.  Hopefully, this will change.  While the concepts might seem orthogonal at first glance, there is a surprising amount of interplay between them.  Beyond Theorem \ref{uniqtrans}, which shows Uniqueness for $\geq \s$ is equivalent to $\lambda$-tameness for basic types, many aspects of frames and frame extensions rely on tameness-like locality principles and, in the other direction, many tameness results, such as categoricity transfer, rely on the concept of minimal types, which were introduced in \cite{sh576} and eventually turned into a frame (see \cite{shelahaecbook}.II.\S3.7).  

It should be noted that there is a loss when these two hypotheses are combined.  We consider here tameness in an AEC with full amalgamation and joint embedding.  These assumptions commonly appear in addition to tameness: amalgamation is used to make types well behaved and joint embedding then follows from $\lambda$-joint embedding.  However, these global assumptions are in contrast to the project of frames, which aims to inductively build up a structure theory, cardinal by cardinal, and derive these properties along the way with the aid of weak diamond.  On the other hand, the existence of frames in the most general setting (see \cite{shelahaecbook}.II.\S3) uses categoricity in two successive cardinals (and more).  If we add no maximal models to this hypothesis, this is  already enough to apply the full categoricity transfer of \cite{tamenessthree}.

On the other hand, the combination of these hypotheses gives much more than just the sum of their parts.  Despite the categoricity transfer results under a tameness hypothesis, there is no robust independence notion for these classes.  The closest approximation is likely Boney and Grossberg \cite{shorttamedep}, where an independence notion of `$< \kappa$ satisfiability' is developed.  Although this notion is well-behaved, additional methods beyond tameness are needed.  Using these method in this paper, we have an independence notion for tame and categorical AECs under some very mild cardinal arithmetic assumptions; see Theorem \ref{forkingincattame}.  Looking at good $\lambda$-frames, the method for building larger frames is a complicated process that changes the Abstract Elementary Class and drops many of the models; see \cite{shelahaecbook}, especially II.\S9.1.  Although this is fine for the end goal, a process that deals with the whole class would likely have more applications.  We provide such a process for tame AECs.

The next section outlines the definitions needed for the rest of the paper, although we assume that the reader is familiar with basic AEC concepts such as Galois types.  Then, the sections show that the various properties of frames extend to $\geq \s$ under the assumption of tameness.  They are organized so that the results only rely on the principles assumed in previous section.  In particular, the stability transfer results of Section \ref{stabsection} do not rely on the tameness for 2-types assumption introduced in Section \ref{symsection}.  We then discuss an application to superstability for AECs in Section \ref{ulmsection} and conclude with an example in Section \ref{hashframesection}.

Important hypotheses are introduced at the end of Sections \ref{prelimsection}, \ref{uniqsection}, and \ref{symsection}.

This paper was written while working on a Ph.D. under the direction of Rami Grossberg at Carnegie Mellon University and I would like to thank Professor Grossberg for his guidance and assistance in my research in general and in this work specifically.  I would also like to thank the referee for their many helpful comments that improved this paper and Alexei Kolesnikov for discussions relating to the example in Section \ref{hashframesection}.
%%%%%%%%%%%%%%%%%%%%%%%%%%%%%%%%%%%%%%%
\section{Preliminaries} \label{prelimsection}
%%%%%%%%%%%%%%%%%%%%%%%%%%%%%%%%%%%%%%%

We assume that the reader is familiar with the definition of AECs and the standard concepts, such as Galois types; see Baldwin \cite{baldwinbook}, Grossberg \cite{grossberg2002}, or \cite{ramibook} for background.  Additionally, frames are covered in depth in Shelah \cite{shelahaecbook}, especially the first volume.  Most of the citations in this paper are from Chapter II of that book, which had previously been circulated as \cite{sh600}.  Occasionally, we will prove a slight variation or weakening of a result from there.  We denote this by adding an asterisk or minus sign, respectively, to the citation and indicate the change.

In order to make this paper more self-contained, we review a few notions in AECs that are not commonplace in the literature.  All of these can be found in the references.
\begin{ndefinition}\label{resdef}
\begin{enumerate}
	\item Given $M \in K$, the set of nonalgebraic types is
	$$\gS^{na}(M) = \{ tp(a/M, N) : tp(a/M, N) \in \gS(M) \te{ and } a \in N - M\}$$
	\item Given $M \in K$ with $\|M\| > LS(K)$ a \emph{resolution} of $M$ is a continuous, strictly increasing sequence of models $\seq{M_i \in K : i < \cf \|M\|}$ so that $M = \cup_{i < \cf \|M\|} M_i$ and $\|M_i\| = |i| + LS(K)$. 
	\item Let $\lambda$ be a cardinal.  $K$ is an \emph{Abstract Elementary Class} in $\lambda$ iff if satisfies every property of being an AEC, \emph{except} we only require it to be closed under chains of length $\leq \lambda$ and all models in $K$ are of size $\lambda$.
\end{enumerate}
\end{ndefinition}

Note that the requirement that all models being the same size makes the axiom about $LS(K)$ meaningless.  The intuition behind this definition is that we sometimes wish to do a local analysis of an AEC by only investigating the models of a particular cardinality; this is the approach that Shelah undertakes with the project of good $\lambda$-frames..  Then $K_\lambda$ is an AEC in $\lambda$.  This technique is useful because we can recover the AEC above $\lambda$ from just this slice at $\lambda$.

\begin{ndefinition}[\cite{shelahaecbook}.II.\S.23] \label{kupdef}
Let $K$ be an AEC in $\lambda$.  We define $(K^{up}, \prec^{up})$ by
\begin{itemize}
	\item $K^{up} = \{ M : M$ is an $L(K)$-structure and there is a directed partial order $I$ and a direct system $\seq{M_s \in K : s \in I}$ such that $M = \cup_{s \in I} M_s\}$
	\item $M \prec^{up} N$ iff there are directed partial orders $I \subset J$ and a direct system $\seq{M_s \in K : s \in J}$ so that $M = \cup_{s \in I} M_s$ and $N = \cup_{s \in J} M_s$.
\end{itemize}
\end{ndefinition}

\begin{nprop}[\cite{shelahaecbook}.II.\S.23]
\begin{enumerate}
	\item If $K$ is an AEC in $\lambda$, then $(K^{up}, \prec^{up})$ is an AEC so $LS(K^{up}) = \lambda$.
	\item If $K$ is an AEC and $\lambda \geq LS(K)$, then $K_\lambda$ is an AEC in $\lambda$ and 
	$$K_{\geq \lambda} = (K_\lambda)^{up}$$
\end{enumerate}
\end{nprop}

We use the definition of frames from Jarden and Shelah \cite{jrsh875}.2.1.1, but with the numbering from Shelah \cite{shelahaecbook}.II.\S2.1.  The missing Axiom (B) is the existence of superlimit model, which is discussed in the Introduction, and Axioms (E)(d) and (i), which are discussed below.

Frames are defined as a triple $\s$ consisting of $K_\lambda$, $\gS^{bs}$, and $\dnf$.  $K_\lambda$ is the collection of all $\lambda$-sized models of some AEC.  $\gS^{bs}$ assigns some well-behaved collection of non-algebraic types to each model in $K_\lambda$ called \emph{basic} types.  $\dnf$ is an independence relation on these basic types and the models of $K_\lambda$.  In the definition and beyond, we will sometimes use the phrase ``$tp(a/M_1, M_3)$ does not fork over $M_0$'' to mean ``$\dnf(M_0, M_1, a, M_3)$ holds.''  This is consistent with the first order terminology and is justified by Axioms (E)(a) and (b).

\begin{ndefinition} \label{framedef}
$\s = (K_\lambda, \dnf_\lambda, \gS^{bs}_\lambda) = (K^\s, \dnf_\s, \gS^{bs}_\s)$ is a good $\lambda$-frame iff
\begin{enumerate}
	
	\item[(A)] $K_\lambda$ is an AEC in $\lambda = \lambda_\s$;
		
	\item[(C)] $K_\lambda$ has AP, JMP, and no maximal models;
	
	\item[(D)] \begin{enumerate}
	
		\item[(a)] $\gS^{bs}_\s(M) \subset \gS(M)$, the domain of $\gS^{bs}_\s$ is $K_\lambda$, and it respects isomorphisms;
		
		\item[(b)] $\gS^{bs}_\s(M) \subset \gS^{na}(M)$;
		
		\item[(c)] {\bf Density:} if $M \precneqq N$ from $K_\lambda$, then there is some $a \in N - M$ so $tp(a/M, N) \in \gS^{bs}_\s(M)$;
		
		\item[(d)] {\bf $bs$-stability:} $|\gS^{bs}_\s(M)| \leq \lambda$ for all $M \in K_\lambda$;
	
	\end{enumerate}

	\item[(E)] \begin{enumerate}
	
		\item[(a)] {\bf Invariance:} $\dnf_\lambda = \dnf_\s = \dnf$ is a four-place relation in which the first, second, and fourth inputs are models from $K_\lambda$ and the third input is an element so $\dnf(M_0, M_1, a, M_3)$ is preserved under isomorphisms and implies i) $M_0 \prec M_1 \prec M_3$; ii) $a \in M_3 - M_1$; and iii) $\dnf(M_0, M_0, a, M_3)$ is equivalent to $tp(a/M_0, M_3) \in S^{bs}_\s(M_0)$;
		
		\item[(b)] {\bf Monotonicity:} if $M_0 \prec M_0' \prec M_1' \prec M_1 \prec M_3'' \prec M_3 \prec M_3'$ and $a \in M_3''$, then $\dnf(M_0, M_1, a, M_3)$ implies $\dnf(M_0, M_1, a, M_3'')$ and $\dnf(M_0', M_1', a, M_3')$;
		
		\item[(c)] {\bf Local Character:} if $\seq{M_i \in K_\lambda : i \leq \delta+1}$ is increasing, continuous, $a \in M_{\delta+1}$, and $tp(a/M_\delta, M_{\delta+1}) \in S^{bs}(M_\delta)$, then there is some $i_0 < \delta$ so $\dnf(M_i, M_\delta, a, M_{\delta+1})$;
		
%		\item {\bf Transitivity:} if $M_0 \prec M_0' \prec M_0'' \prec M_3$ from $K_\alpha$ and $a \in M_3$, then $\dnf(M_0, M_0', a, M_3)$ and $\dnf(M_0', M_0'', a, M_3)$ implies $\dnf(M_0, M_0'', a, M_3)$;

		\item[(e)] {\bf Uniqueness:} If $p, q \in \gS^{bs}(M_1)$ do not fork over $M_0 \prec M_1$ and $p \rest M_0 = q \rest M_1$, then $p = q$;
	
		\item[(f)] {\bf Symmetry:} If $M_0 \prec M_1 \prec M_3$, $a_1 \in M_1$, $tp(a_1/M_0, M_3) \in S^{bs}(M_0)$, and $\dnf(M_0, M_1, a_2, M_3)$, then there are $M_2$ and $M_3'$ so $a_2 \in M_2$, $M_0 \prec M_2 \prec M_3'$, $M_3 \prec M_3'$, and $\dnf(M_0, M_2, a_1, M_3')$;
	
		\item[(g)] {\bf Extension Existence:} If $M \prec N$ and $p \in \gS^{bs}(M)$, then there is some $q \in \gS^{bs}(N)$ so $p \leq q$ and $q$ does not fork over $M$;

		\item[(h)] {\bf Continuity:} if $\seq{M_i \in K_\lambda : i \leq \delta}$ with $\delta$ limit, $p \in M_\delta$, and, for all $i < \delta$, $p \rest M_i$ does not fork over $M_0$, then $p \in \gS^{bs}(M_\delta)$ and $p$ does not fork over $M_0$.
		
%		\item {\bf Non-forking Amalgamation:} if, for $\ell = 1, 2$, $M_0 \prec M_\ell$ from $K_\lambda$, $a_\ell \in M_\ell - M_0$, and $tp(a_\ell/M_0, M_\ell) \in S^{bs}(M_0)$, then there are $f_1, f_2, M_3$ so $M_0 \prec M_3 \in K_\lambda$ and, for $\ell = 1, 2$, we have $f_\ell:M_\ell \to_{M_0} M_3$ and $\dnf(M_0, f_{3-\ell}(M_{3-\ell}), f_\ell(a_\ell), M_3$.

	\end{enumerate}

\end{enumerate}

\end{ndefinition}

Note that basic types are types of singletons and so all parameters are single elements.  Thus, references to basic types always implicitly refer to 1-types.  In \cite{shelahaecbook}.II.\S2, Shelah points out that the definition could be changed to finite types ``with no real loss;'' in this case, the results of this paper for extending frames to higher cardinals would require tameness for $<\omega$-types.  Shelah makes this comment explicit in \cite{shelahaecbook}.III.\S5.2 with the notion of independence sets, which gives a notion of nonforking for larger and possibly infinite tuples.  This is further studied by Jarden and Sitton in \cite{jardensitton}.  Using this definition, the crucial 2-types in Theorem \ref{symtheorem}--$tp(a_1 a'/M[0, 0, \mu], M[2, 1, \mu])$ in the notation there--would be basic 2-types.  This will be explored further in \cite{extendingtameness}.

Typically, when we cite the frame axioms, we will do so by letter in theorem statements and by name elsewhere.  Also, the unnamed axioms (Axioms (A) and (D)(a) and (b)) and Invariance are taken to be basic, so we will often not mention them even from lists of axioms used in a proof.  This is because they are satisfied of all possible candidates for independence relations.

Several examples of frames are given in \cite{shelahaecbook}.II.\S3.  For the first order case, forking in superstable theories satisfies the definition with regular types.  In a general AEC, Shelah \cite{sh576} derives a good $\lambda$-frame from categoricity in two successive cardinals, a medium number in the next, and some cardinal arithmetic; see \cite{shelahaecbook}.VII.\S0.3, .4 for definitions or Shelah \cite{shelahcardinalarithmetic} for a larger discussion.

\begin{ntheorem}[\cite{shelahaecbook}.II.\S3.7] \label{genframeext}
Assume $2^\lambda < 2^{\lambda^+} < 2^{\lambda^{++}}$ and
\begin{enumerate}
	\item $K$ is an AEC with $LS(K) \leq \lambda$;
	\item $K$ is categorical in $\lambda$ and $\lambda^+$;
	\item $K$ has a model in $\lambda^{++}$; and
	\item $I(\lambda^{++}, K) < \mu_{unif}(\lambda^{++}, 2^{\lambda^+})$ and $WDmId(\lambda^+)$ is not $\lambda^{++}$-saturated.
\end{enumerate}
Then there is a good $\lambda^+$-frame.
\end{ntheorem}

Shelah shows the following additional properties hold of frames, which he originally states as axioms.

\begin{ntheorem}[\cite{shelahaecbook}.II.\S2.18, .16]
\begin{itemize}

	\item Axioms (A), (C), (D)(a) and (b), and (E)(a), (b), (e), and (g) imply Axiom (E)
	\begin{enumerate}
	\item[(d)] {\bf Transitivity:} if $M_0 \prec M_0' \prec M_0'' \prec M_3$ from $K_\alpha$ and $a \in M_3$, then $\dnf(M_0, M_0', a, M_3)$ and $\dnf(M_0', M_0'', a, M_3)$ implies $\dnf(M_0, M_0'', a, M_3)$.
	\end{enumerate}

	\item Axioms (A), (C), and (E)(b), (d), (f), and (g) imply Axiom (E)
	\begin{enumerate}
	\item[(i)] {\bf Non-forking Amalgamation:} if, for $\ell = 1, 2$, $M_0 \prec M_\ell$ from $K_\lambda$, $a_\ell \in M_\ell - M_0$, and $tp(a_\ell/M_0, M_\ell) \in \gS^{bs}(M_0)$, then there are $f_1, f_2, M_3$ so $M_0 \prec M_3 \in K_\lambda$ and, for $\ell = 1, 2$, we have $f_\ell:M_\ell \to_{M_0} M_3$ and $\dnf(M_0, f_{3-\ell}(M_{3-\ell}), f_\ell(a_\ell), M_3)$.
	\end{enumerate}
	
\end{itemize}
\end{ntheorem}

We conclude by Shelah's exercise in increasing the size of frames.  This can be seen as a generalization of the standard technique of taking an AEC in $\lambda$ and blowing it up to an AEC; see \cite{shelahaecbook}.II.\S1.23.  We replace his notation ``$\dnf_{< \infty}$'' with ``$\dnf_{\geq \s}$'' because it is more consistent with the notion of referring to the extended frame $\geq \s = (K, \gS^{bs}_{\geq \s}, \dnf_{\geq \s})$.

\begin{ndefinition} \label{upframes}
\begin{enumerate}

	\item[]
	
	\item[\cite{shelahaecbook}.II.\S2.4.1)] $K^{3, bs} = K_{\geq \mathfrak{s}}^{3, bs} = \{ (a, M, N) \in K^{3, na}  :$there is $M' \prec M$ from $K_\lambda$ such that, for all $M'' \in K_\lambda$, $M' \prec M'' \prec M$ implies that $tp(a/M'', N) \in \gS^{bs}(M'')$ does not fork over $M' \}$

	\item[\cite{shelahaecbook}.II.\S2.7/.8.1)] For $M \in K$,
	$$\gS_{\geq s}^{bs}(M) = \{ p \in \gS(M) : \te{for some/every $tp(a/M, N) = p$, $(a, M, N) \in K^{3, bs}$ }\}$$
	
	\item[\cite{shelahaecbook}.II.\S2.5)] We say that $\dnf_{\geq \s}(M_0, M_1, a, M_3)$ holds iff $M_0 \prec M_1 \prec M_3 \in K$, $a \in M_3 - M_1$, and there is $M_0' \prec M_0$ from $K_\lambda$ such that if $M_0' \prec M_1' \prec M_1$ and $M_1' \cup \{a\} \subset M_3' \prec M_3$ with $M_1', M_3' \in K_\lambda$, then $\dnf_{\s}(M_0', M_1', a, M_3')$.

	\item[\cite{shelahaecbook}.II.\S2)] If $\s$ is a good $\lambda$-frame, then set $\geq \s := ((K_\s)^{up}, \gS^{bs}_{\geq \s}, \dnf_{\geq \s})$.
	
	\item $\geq \s$ is a \emph{good frame} iff it satisfies the axioms for good $\lambda$-frames after removing the restriction on the size of the models and length of sequences.

\end{enumerate}
\end{ndefinition}

Many of the properties of good $\lambda$-frames transfer upwards immediately.

\begin{ntheorem} \label{basicextension}
If $\s$ is a good $\lambda$-frame, then $\geq \s$ is a good frame, except possibly for (C), (D)(d), and (E)(e), (f), and (g).
\end{ntheorem}

{\bf Proof:} By the results of \cite{shelahaecbook}.II.\S2.  Specifically, Invariance and (D)(a) are 8.3, Density is 9, Monotonicity is 11.3, Transitivity is 11.4, Local Character is 11.5, and Continuity is 11.6.\hfill \dag\\

At least some additional hypothesis is necessary to transfer all properties of a good $\lambda$-frame $\s$ to a good frame $\geq \s$.  This can be observed by observing that the Hart-Shelah examples \cite{hash323} (reanalyzed more deeply by Baldwin and Kolesnikov in \cite{untame}) have good $\lambda$-frames at low cardinalities, but the upward extension fails Uniqueness and Basic Stability (and only those) exactly at the cardinal that tameness breaks down; see Section \ref{framehashsection} for details.

In light of this, to prove that $\geq \s$ is a good frame, we need to additionally show amalgamation, joint embedding, no maximal models, uniqueness, basic stability, extension existence, and symmetry.  In order to avoid any mention of categoricity or non-structure arguments that require instances of the weak continuum hypothesis (as in \cite{sh576} or \cite{shelahaecbook}.I.\S3), we assume amalgamation.  This leads us to our first hypothesis.

\begin{hyp} \label{basichyp}
$K$ is an AEC with $LS(K) \leq \lambda = \lambda_\s$ with amalgamation and $\s$ is a good $\lambda$-frame.
\end{hyp}

Although joint embedding is not included in this hypothesis, we may freely use it due to the following fact.

\begin{nfact}
If $K$ is an AEC with amalgamation and $K_\lambda$ has joint embedding, then $K_{\geq \lambda}$ has joint embedding.
\end{nfact}

Additionally, Jarden and Shelah \cite{jrsh875} introduce the notion of semi-good $\lambda_\s$-frames, which replace Basic Stability with Almost Basic Stability, which requires that $|\gS^{bs}_\s(M)| \leq \lambda_\s^+$ for all $M \in K_{\lambda_\s}$.  The following could also be done for semi-good frames, although Section \ref{stabsection} shows that, assuming tameness, $\geq \s$ will be stable everywhere strictly above $\lambda_\s$, even if $\s$ is just a semi-good $\lambda_\s$-frame.

%%%%%%%%%%%%%%%%%%%%%%%%%%%%%%%%%%%%%%%
\section{Tameness and Uniqueness} \label{uniqsection}
%%%%%%%%%%%%%%%%%%%%%%%%%%%%%%%%%%%%%%%

Tameness is the key property that is necessary in extending frames, needed both for Uniqueness and Symmetry.  In this section, we show that tameness for 1 types is \emph{equivalent} to the frame having uniqueness.  Recall the definition of tameness.

\begin{ndefinition}
We say that $K$ is $(\lambda, \kappa)$ tame for $\alpha$-types iff, given any $M \in K_\kappa$ and $p \neq q \in \gS(M)$ of length $\alpha$, there is some $N \prec M$ of size $\lambda$ so $p \rest N \neq q \rest N$.

We say that $K$ is $\lambda$ tame for $\alpha$-types iff it is $(\lambda, \kappa)$ tame for all $\kappa \geq \lambda$.

If we omit the $\alpha$, then we mean 1-types.

Tameness for basic types is the same property with $p \neq q \in \gS_{bs}(M)$.
\end{ndefinition}

We will use this only for $\alpha$ equal to 1 (this section) or 2 (Section \ref{symsection}).

We will prove the following.

\begin{ntheorem} \label{uniqtrans}
$K_{\geq \mathfrak{s}}$ is $\lambda_\s$-tame for basic types iff $\geq \mathfrak{s}$ satisfies Uniqueness.
\end{ntheorem}

We can parametrize this result and get that $(\lambda_\s, \mu)$-tameness is equivalent to Uniqueness for models of size $\mu$.  To prove this, we use and prove the following variation of a claim from Shelah's book:

\begin{nclaim}[\cite{shelahaecbook}.II.\S2.10*] \label{2.10}
If $tp(a/M, N) \in \gS^{\te{bs}}_{\geq \mathfrak{s}}(M)$, then there is some $M_0 \in K_\s$ with $M_0 \prec M$ such that $tp(a/M_0, N) \in \gS_\mathfrak{s}^{\te{bs}}(M_0)$ and
\begin{eqnarray*}\te{ if $M_0 \prec M' \prec M$, then $tp(a/M', N) \in \gS_{\geq \mathfrak{s}}^{\te{bs}}(M')$ does not fork over $M_0$. }\end{eqnarray*}
\end{nclaim}

Our trivial variation is to permit $M' \in K_{\geq \s}$, rather than restricting to the case $M' \in K_\s$.  We offer a proof for completeness because \cite{shelahaecbook} omits one.

{\bf Proof:} Since $tp(a/M, N) \in \gS_{\geq \s}^\bs(M)$, we have that $(a, M, N) \in K^{3, \bs}_{\geq s}$.  So, by definition, there is some $M_0 \prec M$ of size $\lambda = \lambda_\s$ so, for all $M'' \in K_\lambda$,
$$M_0 \prec M'' \prec M \implies \dnf_{\s}(M_0, M'', a, N)$$
Now we just need to prove (1).\\
Let $M' \in K$ such that $M_0 \prec M' \prec M$.  First, we see that $(a, M', N) \in K^{3, \bs}_{\geq \s}$ as witnessed by $M_0$.  Now we want to show that $\dnf_{\geq \s} (M_0, M', a, N)$ and, in fact, we claim that $M_0$ is the witness $M_0'$ for this.  If $M_1' \in K_\lambda$ such that $M_0 \prec M_1' \prec M'$, then, since $M' \prec M$, $M_1' \prec M$.  Thus, by the definition of $M_0$ as the witness for $(a, M, N) \in K_{\geq \s}^{3, \bs}$, $tp(a/M_1', N) \in \gS_\s^\bs(M_1')$ does not fork over $M_0$.  So then $\dnf_{\geq \s} (M_0, M', a, N)$ as desired. \hfill \dag\\

We define an equivalence relation $\mathcal{E}^\s_M$, as in \cite{shelahaecbook}.II.\S2.7.3, for $M \in K_{\geq \s}$ on $\gS^\bs_{\geq \s}(M)$ by $p \mathcal{E}^\s_M q$ iff $p \rest N = q \rest N$ for all $N \prec M$ in $K_\s = K_{\lambda_\s}$.

We quote:

\begin{nfact}[\cite{shelahaecbook}.II.\S2.8.5]
$\mathcal{E}^\s_M$ is an equivalence relation on $\gS_{\geq \s}^\bs(M)$ and if $p, q \in \gS_{\geq \s}^\bs(M)$ do not fork over $N \in K_\s$ such that $N \prec M$ then
$$p \mathcal{E}^\s_M q \iff (p \rest N = q \rest N)$$
\end{nfact}

{\bf Proof of Theorem \ref{uniqtrans}:}  First, suppose that $\geq \s$ satisfies Uniqueness for some $M \in K_\mu$ with $\mu \geq \lambda_\s$.  Let $p, q \in \gS_{\geq \s}^\bs(M)$ such that $p \rest N = q \rest N$ for all $N \prec M$ of size $\lambda$.  Then we can find $M_p, M_q$ as in Claim \ref{2.10} above.  Let $M' \prec M$ of size $\lambda$ contain both.  Then by Monotonicity, we know that $p$ and $q$ both don't fork over $M'$.  However, by assumption, $p \rest M' = q \rest M'$.  Then, by Uniqueness, $p = q$.\\
Second, suppose that $K_s$ is $(\lambda_\s, \mu)$ tame for basic types.  In particular, this means that $\mathcal{E}^\s_M$ is equality for all $M \in K_\mu$.  Let $M \in K_\mu$, $p, q \in \gS_{\geq \s}^\bs(M)$, and $M' \prec M$ such that $p$ and $q$ do not fork over $M'$ (in the sense of $\geq \s$) and $q \rest M' = q \rest M'$.  By Claim \ref{2.10}, there are $M_p, M_q \prec M$ of size $\lambda$ such that $p\rest M'$ does not fork over $M_p$ and $q \rest M'$ does not fork over $M_q$.  As above, find $M_0 \prec M'$ of size $\lambda$ to contain $M_p$ and $M_q$; then by Monotonicity, $p \rest M'$ and $q \rest M'$ do not fork over $M_0$.  Then by Transitivity, $p$ and $q$ don't fork over $M_0$.  Also, since $M_0 \prec M'$ and $p \rest M' = q \rest M'$, we have $p \rest M_0 = q \rest M_0$.  Then, by \cite{shelahaecbook}.II.\S2.8.5, $p \mathcal{E}^\s_M q$.  But, by tameness, this is equality, so $p = q$. \hfill \dag\\

Additionally, if $\s$ is type full ($S^{bs}_\s = S^{na}$), then we can extend our result on tameness to not mentioning basic types at all.

\begin{ncor}
If $\geq \s$ is a type full good frame, then $K_{\geq \s}$ is $\lambda_\s$ tame.
\end{ncor}

Note that \cite{shelahaecbook}.II.\S6.36 says that we can assume $\s$ is full if it has existence for $K^{3, uq}_{\lambda_\s}$ (see \cite{shelahaecbook}.II.\S5.3).

In light of these results, we add the following hypothesis.

\begin{hyp} \label{tamehyp}
$K$ is $\lambda_\s$-tame for basic 1-types.
\end{hyp}

%%%%%%%%%%%%%%%%%%%%%%%%%%%%%%%%%%%%%%%
\section{Basic Stability} \label{stabsection}
%%%%%%%%%%%%%%%%%%%%%%%%%%%%%%%%%%%%%%%

In this section, we use only tameness for 1-types (and therefore no symmetry) to prove that an extended frame leads to basic stability in all larger cardinals.  This is similar to the first order argument that stability and $\kappa(T) = \omega$ together imply superstability.  This has been done in non-elementary contexts by Makkai and Shelah \cite{makkaishelah}.4.14.

\begin{ntheorem}\label{stabtransthm}
For all $\kappa \geq \lambda$, $K$ is $\kappa$-stable for $\geq \s$ basic types; that is, for all $M \in K_\kappa$, $|\gS^{bs}_{\geq \s}(M)| \leq \kappa$.
In particular, $(\lambda, \leq \kappa)$-tameness for basic 1-types implies $\kappa$-stability for basic types. 
\end{ntheorem}

{\bf Proof:} We proceed by induction on $\lambda \leq \mu \leq \kappa$.  If $\mu = \lambda$, then this is the hypothesis.  For $\mu > \lambda$, let $M \in K_\mu$ and find a resolution $\seq{M_i \in K_{[\lambda, \mu)} : i < \cf \mu}$ of $M$.  By Local Character for $\geq \s$, for each $p \in S^{bs}_{\geq \s} (M)$, there is some $i_p < \cf \mu$ such that $p$ does not fork over $M_{i_p}$ in the sense of $\geq \s$.  By Theorem \ref{uniqtrans}, $\geq \s$ satisfies Uniqueness for domains of size $\mu$, so the map $p \mapsto M_{i_p}$ is injective from $\gS^{bs}_{\geq \s}(M)$ to $\bigcup_{i < \cf \mu} \gS^{bs}_{\geq \s}(M)$. So
$$|\gS^{bs}_{\geq \s}(M)| \leq \sum_{i < \cf \mu} |\gS^{bs}_{\geq \s}(M_i)| = \mu$$
\hfill \dag\\

In particular, this uses only Local Character and Uniqueness.  We can extend this to full stability using \cite{shelahaecbook}.II.\S.4.2.1, which shows that stability for basic types implies stability for all types using amalgamation, Density, Monotonicity, and Local Character.  Thus, we get the following stability transfer which improves on results of Grossberg and VanDieren \cite{tamenessone}; Baldwin, Kueker, and VanDieren \cite{bkv}; and Lieberman \cite{liebermanrank}, but adds the assumption of a good $\lambda$-frame.

\begin{ncor} \label{stabtrans}
Suppose $K$ is $\chi$-tame for 1-types and has a good $\chi$-frame except possibly for the assumption of basic stability.  If $K$ is stable (or just stable for basic types) in some $\kappa \geq \chi$, then it is stable in all $\mu \geq \kappa$.
\end{ncor}

%%%%%%%%%%%%%%%%%%%%%%%%%%%%%%%%%%%%%%%
\section{Extension Existence} \label{extsection}
%%%%%%%%%%%%%%%%%%%%%%%%%%%%%%%%%%%%%%%

We now turn to the existence of nonforking extensions of basic types.  One of the difficulties of using Galois types (compared to syntactic types) is that an increasing sequence of types need not have an upper bound.  Shelah and Baldwin \cite{nonlocality}.3.3 construct an example of an AEC that has an increasing sequence of types with no upper bound from $2^{\aleph_0} = \aleph_1$, $\diamondsuit_{\aleph_1}$, $\diamondsuit_{S^{\aleph_2}_{\cf \aleph_1}}$, and $\square_{\aleph_2}$.  However, if we require that the sequence is \emph{coherent} (see below), then there is an upper bound.  Equivalently, Shelah \cite{sh576} and others work with increasing sequences from $K^{3, na}_\lambda$.

In essence, we will show that a good $\lambda$-frame and $\lambda$-tameness imply that types are local and apply an argument similar to \cite{sh394} (proved as \cite{baldwinbook}.11.5) to show that compactness follows; see \cite{nonlocality} for the relevant definitions, although we will not use them here.  This is essentially the same argument used in the proof of \cite{tamenessthree}.2.22, where they work with quasiminimal types.  In all cases, there is some property--locality, quasiminimality, or Uniqueness--that ensures that there is only one possible extension at limit steps.  We reprove this here because previous contexts have worked inside of a monster model, which we do not have.  However, the ideas in Proposition \ref{framecoher} are not new.

\begin{ndefinition}
Given increasing sequences $\seq{M_i : i < \delta}$ and $\seq{p_i \in \gS(M_i) : i < \delta}$, the sequence of types is called \emph{coherent} iff there are, for $j < i < \delta$, models $N_i$, elements $a_i$, and maps $f_{j, i}:N_j \to N_i$ so
\begin{enumerate}

	\item for all $k < j < i < \delta$, we have $f_{k, i} = f_{j, i} \circ f_{k, j}$;
	
	\item $tp(a_i/M_i, N_i) = p_i$;
	
	\item $f_{j, i}$ fixes $M_j$ for all $i > j$; and
	
	\item $f_{j, i}(a_j) = a_i$.

\end{enumerate}
\[
\xymatrix{
(N_k, a_k) \ar@/^2pc/[rr]^{f_{k, i}} \ar[r]^{f_{k, j}} & (N_j, a_j) \ar[r]^{f_{j, i}} & (N_i, a_i)\\
M_k \ar[r] \ar[u] & M_j \ar[r] \ar[u] & M_i \ar[u]
}
\]

\end{ndefinition}

If we have a coherent sequence of types, it must have an upper bound.  Namely, taking $M = \bigcup_{i < \delta} M_i$, $(N, f_i^*)_{i < \delta} = \varinjlim_{j < k < \delta} (N_k, f_{j, k})$, and $a^* = f^*_0(a_0)$, the upper bound is $tp(a^*/M, N)$.

The above does not require frames.  However, if we have a frame, then all nonforking sequences of types are coherent.

\begin{nprop} \label{framecoher}
Let $\seq{M_i \in K_{\geq \lambda_\s} : i < \delta}$ be an increasing, continuous sequence.  If $\seq{p_i \in \gS^{bs}_{\geq \s}(M_i) : i < \delta}$ is an increasing sequence of basic 1-types such that each $p_i$ does not fork over $M_0$, then $p_i$ is coherent.  Thus, there is $p_\delta \in S^{bs}_{\geq \delta}(M_\delta)$ extending each $p_i$.
\end{nprop}

Note that Uniqueness (which follows from Theorem \ref{uniqtrans} since we assumed tameness for basic types in Hypothesis \ref{tamehyp}) is the key property used in this proof.

{\bf Proof:}  For $i = 0$, set $(a_0, M_0, N_0) \in K^{3, bs}$ to be some triple realizing $p_0$.

For $i$ limit, set $(N_i, f_{j, i})_{j < i} = \varinjlim_{l < k < i}(N_k, f_{l, k})$.  Then $M_i \prec N_i$.  Set $a_i = f_{0, i}(a_i)$, which is equal to $f_{j, i}(a_j)$ for any $j < i$.  For each $j < i$, $f_{j, i}$ fixes $M_j$, so $a_i \vDash p_j$.  Thus, $tp(a_i/M_j, N_i)$ doesn't fork over $M_0$.  Since this is true for all $j < i$, Continuity says that $tp(a_i/M_i, N_i)$ does not fork over $M_0$.  Since $p_i$ also does not fork over $M_0$, Uniqueness implies that $tp(a_i/M_i, N_i) = p_i$, as desired.

For $i = j+1$, find $(a_i', M_i, N_i')$ such that $tp(a_i'/M_i, N_i') = p_i$.  Since $p_i \rest M_j = p_j$, $a_i$ and $a_j$ realize the same type over $M_i$.  Thus we can construct the following commutative diagram

\[
 \xymatrix{N_j  \ar[rr]^{f_{j, i}}  & & N_i\\
M_j \ar[u] \ar[r] & M_i \ar[r] & N_i' \ar[u]
 }
\]

so $f_{j, i}(a_j) = a_i$.  Then set $f_{k, i} = f_{j, i} \circ f_{k, j}$ for any $k < j$. 

Once we have constructed the coherent sequence, there is some $p \in \gS(M)$ for $M = \cup_{i < \delta} M_i$ that extends each $p_i$.  By Continuity, $p \in \gS^{bs}_{\geq \s}(M)$ and $p$ does not fork over $M_0$.\hfill \dag\\

Now we prove that Extension Existence holds in $\geq \s$.  We proceed by induction.

\begin{ntheorem} \label{exttrans}
$\geq \s$ satisfies Axiom (E)(g).
\end{ntheorem}

{\bf Proof:} We want to show:

\begin{center}
If $M \prec N$ from $K_{\geq \lambda_\s}$ and $p \in \gS_{\geq \s}^\bs(M)$, then there is some $q \in \gS_{\geq \s}^\bs(N)$ such that $p \leq q$ and $q$ does not fork over $M$ (in the $\dnf_{\geq \s}$ sense).
\end{center}

We will prove this by induction on $\|N\|$.\\

{\bf Base Case:} $\| N \| = \lambda_\s$\\
Then $\| M \| = \lambda_\s$ as well, and this follows from $\s = (\geq \s) \rest \lambda_\s$ being a good $\lambda_\s$-frame.
	
{\bf Inductive Step:} $\| N \| = \mu > \lambda_\s$\\
We break into two cases based on the size of $M$.\\
If $\|M\|<\|N\|$, then we find a resolution $\seq{N_i \in K_{< \mu} \mid i < \mu}$ such that $N_0 = M$.  By induction, we will construct increasing $p_i \in \gS_{\geq \s}^\bs(N_i)$ such that $p_{i}$ does not fork over $N_0$ and extends $p$.  Clearly, $p_0 = p$.\\
		For $i$ limit, by Proposition \ref{framecoher}, we can find some $p_i$ such that $p_i \rest N_j = p_j$ for all $j < i$.  Then $p_i \rest N_j$ does not fork over $M$ for all $j < i$, so, by Continuity, $p_i$ does not fork over $M$.\\
		For $i = j + 1$, we use our induction to extend $p_j$ to some $p_i \in \gS_{\geq \s}^\bs(N_i)$ that doesn't fork over $N_j$; this is valid since $\| N_i \| < \| N \|$.\\
		Then, we use Proposition \ref{framecoher} a final time to find $q \in \gS_{\geq \s}^\bs(N)$ such that $q \rest N_i = p_i$.  By Continuity, this means $q$ does not fork over $M$ as desired.\\
		If $\|M\| = \|N\|$, we find $M_0 \prec M$ in $K_\s$ as in Claim \ref{2.10} such that if $M_0 \prec M' \prec M$, $p \rest M'$ does not fork over $M_0$.  Then we use this as the start for a resolution $\seq{M_i \in K_{< \mu} \mid i < \cf \mu}$ of $M$.  Set $p_i = p \rest M_i$; note that $p_i$ does not fork over $M_0$.  Now we find a resolution $\seq{N_i \in K_{< \mu} \mid i < \cf \mu}$ of $N$ such that $M_i \prec N_i$.  We are going to find increasing $q_i \in \gS_{\geq \s}^\bs(N_i)$ by induction such that $q_i$ does not fork over $M_0$ and $p_i \leq q_i$.\\
			We use the induction hypothesis to find $q_0 \in \gS_{\geq \s}^\bs(N)$ that extends $p_0$ and does not fork over $M_0$.\\
			For $i$ limit, use the induction hypothesis to find $q_i \in \gS_{\geq \s}^\bs(N_i)$ that extends all $q_i$.  By continuity, $q_i$ does not fork over $M$ or over $N_j$ for all $j < i$.\\
			For $i = j + 1$, use induction to find $q_i \in \gS_{\geq \s}^\bs(N_i)$ such that $q_i \geq q_j$ and $q_i$ does not fork over $N_j$.  Then, by Transitivity, $q_i$ does not fork over $M_0$.  Also note that $p_i$ does not fork over $M_0$ and $q_i \rest M_0 = p_0 = p_i \rest M_0$, so Uniqueness tells us that $q_i \rest M_i = p_i \rest M_i = p_i$.\\
			Now we use Proposition \ref{framecoher} to set $q \in \gS_{\geq \s}^\bs(N)$ to extend all $q_i$ and $p_0$.  Again by Continuity, $q$ does not fork over $M_0$.  Also, $q\rest M_0 = p_0 = p \rest M_0$ so, since $p$ also does not fork over $M_0$, we can use Uniqueness to get that $q \rest M = p$.  Finally, by Monotonicity, we have that $q$ does not fork over $M$. \hfill \dag\\			
	
%%%%%%%%%%%%%%%%%%%%%%%%%%%%%%%%%%%%%%%
\section{Tameness and Symmetry} \label{symsection}
%%%%%%%%%%%%%%%%%%%%%%%%%%%%%%%%%%%%%%%

In this section, we show that tameness for 2-types implies Symmetry in $\geq \s$.  Unfortunately, unlike Section \ref{uniqsection}, this is not shown to be an equivalence.  This is enough for our goal of extending a frame, but a characterization of exactly when Symmetry holds in $\geq \s$ would be better.  We know that tameness for 2-types (or even tameness for basic 2-types in the sense of \cite{shelahaecbook}.III.\S5.2) is not this characterization because the Hart-Shelah examples have frames with Symmetry at all cardinals, including after the tameness fails; see Section \ref{hashframesection}.  Additionally, the precise relationship between tameness for 1-types and tameness for 2-types is not currently known, although tameness for 2-types clearly implies tameness for 1-types.

\begin{ntheorem} \label{symtheorem}
If $K$ satisfies $\lambda_\s$ tameness for 2-types, then $\geq \s$ satisfies Axiom (E)(f).
\end{ntheorem}

For reference, a diagram of the models involved in the proof is included in the proof.  This diagram and the naming convention for models deserves some explanation and we are indebted to the referee for pushing us to a better presentation.  Functions like $f$ and $g$ \emph{above} arrows have their usual meanings (that $f$ is a $K$ embedding between models), but we write elements \emph{under} arrows to indicate that the element is in the end model but not the starting model.  

The majority of models are names $M[i, j, \chi]$ for $i$ and $j$ natural numbers and $\chi$ a cardinal, either $\mu$ or $\lambda_\s$.  The cardinal $\chi$ denotes the size of the model and the sizes separate the models into two levels.  If we have $M[i, j, \chi]$ and $M[i', j', \chi']$ with $i \leq i'$, $j \leq j'$, and $\chi \leq \chi'$, then this will mean that $M[i, j, \chi]$ is embedded into $M[i', j', \chi']$ by the map indicated by the diagram.  In particular, we \emph{do not} have an embedding of $M[4, 2, \lambda_\s]$ into $M[3, 1, \mu]$, even though the first model is below the second in the diagram (and the nonstandard indices for the first model are picked to indicate this).  The exception to this convention are the models $M^-$, $M_I$, and $M_{II}$.  These models are all contained in $M[0, 0, \lambda_\s]$ and are used as ``helper models;'' that is, they lend properties to $M[0, 0, \lambda_\s]$, but are not directly used in the proof.

{\bf Proof:}  Suppose we have $M[0, 0, \mu], M[0, 1, \mu], M[1, 1, \mu] \in K_\mu$ such that  $\dnfty(M[0,0, \mu], M[0, 1, \mu], a_1, M[1, 1, \mu])$ and $a_2 \in M[0, 1, \mu]$ such that $tp(a_2/M[0,0, \mu], M[1, 1, \mu]) \in \gS^{bs}_{\geq \s}(M[0,0, \mu])$.  Let $M[1, 0, \mu] \in K_\mu$ such that $M[0,0, \mu] \prec M[1, 0, \mu] \prec M[1, 1, \mu]$ and $a_1 \in M[1, 0, \mu]$.  By Extension Existence, there is $M[2, 1, \mu] \in K_\mu$ such that $M[1, 1, \mu] \prec M[2, 1, \mu]$ and $a' \in M[2, 1, \mu]$ such that $\dnfty(M[0,0, \mu], M[1, 0, \mu], a', M[2, 1, \mu])$ and $tp(a'/M[0,0, \mu], M[2, 1, \mu]) = tp(a_2/M[0,0, \mu], M[1, 1, \mu])$.  We want to show that this type equality still holds if we add $a_1$.\\
{\bf Main Claim:} $tp(a_1 a_2/M[0,0, \mu], M[1, 1, \mu]) = tp(a_1 a'/M[0,0, \mu], M[2, 1, \mu])$\\
{\bf This is Enough:} Let $N \in K_\mu$ witness the above type equality; that is, $M[1, 1, \mu] \prec N$ and there is $f:M[2, 1, \mu] \to_{M[0,0, \mu]} N$ such that $f(a_1 a') = a_1 a_2$.  Then apply $f$ to $\dnfty(M[0,0, \mu], M[1, 0, \mu], a', M[2, 1, \mu])$; this shows that $\dnfty(M[0,0, \mu], f(M[1, 0, \mu]), a_2, N)$.  This proves Symmetry since $a_1 \in f(M[1, 0, \mu])$.\\

{\bf Proof of Main Claim:} Fix $M^- \prec M[0,0, \mu]$ of size $\lambda_\s$.  From the assumption of tameness for 2-types, it suffices to show 
$$tp(a_1a_2/M^-/M[1, 1, \mu]) = tp(a_1 a'/M^-, M[2, 1, \mu])$$
By the definition of $\geq \s$, there are $M_I, M_{II} \in K_{\lambda_\s}$ such that $M_I, M_{II} \prec M[0,0, \mu]$ and that witness (in the sense of the definition of $\dnf_{\geq \s}$, see Definition \ref{upframes}) $\dnfty(M[0,0, \mu], M[0, 1, \mu], a_1, M[1, 1, \mu])$ and $\dnfty(M[0,0, \mu], M[1, 0, \mu], a', M[2, 1, \mu])$, respectively.  Let $M[0,0,\lambda_\s]\in K_{\lambda_\s}$ such that $M[0, 0, 0]_{\lambda_\s} \prec M[0,0, \mu]$ and it contains $M^-, M_I$, and $M_{II}$.  Then, since $M[0, 0, \lambda_\s]$ contains witnesses to the nonforking, we have that
\begin{enumerate}

	\item if there are $M, M' \in K_{\lambda_\s}$ with $a_1 \in  M'$ such that $M[0,0,\lambda_\s] \prec M \prec M[0, 1, \mu]$ and $M \prec  M' \prec M[1, 1, \mu]$, then $\dnf_{\s}(M[0,0,\lambda_\s], M, a_1, M')$; and
	
	\item if there are $M, M' \in K_{\lambda_\s}$ with $a' \in M'$ such that $M[0,0,\lambda_\s] \prec M \prec M[1, 0, \mu]$ and $ M \prec  M' \prec M[2, 1, \mu]$ , then $\dnf_{\s}(M[0,0,\lambda_\s], M, a',  M')$.

\end{enumerate}

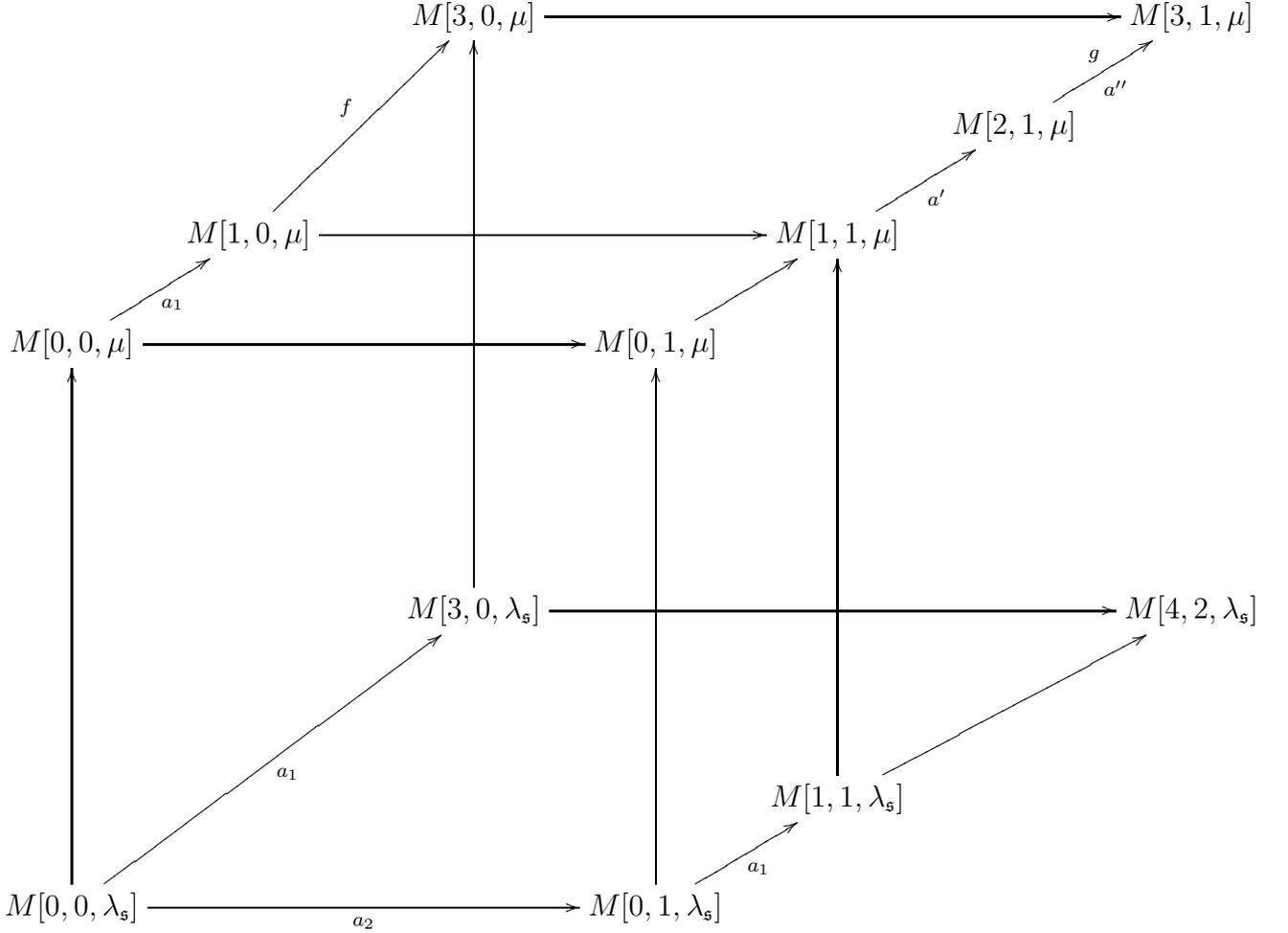
\[
\xymatrix@C=1em{
 &  &  & M[3, 0, \mu] \ar[rrrr] &  &  &  & M[3, 1, \mu] \\
 &  &  &  &  &  & M[2, 1, \mu] \ar[ur]_{a''}^{g} &  \\
 & M[1, 0, \mu] \ar[rrrr] \ar[uurr]^{f} &  &  &  & M[1, 1, \mu] \ar[ur]_{a'} &  &  \\
M[0, 0, \mu] \ar[ur]_{a_1} \ar[rrrr] &  &  &  & M[0, 1, \mu] \ar[ur] &  &  &  \\
 &  &  &  &  &  &  &  \\
 &  &  &  &  &  &  &  \\
 &  &  & M[3, 0, \lambda_\s] \ar[uuuuuu] \ar[rrrr] &  &  &  & M[4, 2, \lambda_\s] \\
 &  &  &  &  &  &  &  \\
 &  &  &  &  & M[1, 1, \lambda_\s] \ar[uurr] \ar[uuuuuu] &  &  \\
M[0, 0, \lambda_\s] \ar[uuuuuu] \ar[rrrr]_{a_2} \ar[uuurrr]_{a_1} &  &  &  & M[0, 1, \lambda_\s] \ar[uuuuuu] \ar[ur]_{a_1} &  &  &  \\
}
\]

Since $\lambda_\s \geq LS(K)$, there are $M[0, 1, \lambda_\s], M_[1, 1, \lambda_\s] \in K_{\lambda_\s}$ such that $M[0,0,\lambda_\s] \prec M[0,1,\lambda_\s] \prec M[0, 1, \mu]$ and $a_2 \in M[0, 1, \lambda_\s]$; and $M[0,1,\lambda_\s] \prec M[1, 1, \lambda_\s] \prec M[1, 1, \mu]$ and $a_1 \in M[1, 1, \lambda_\s]$.  From the definition of $M[0,0,\lambda_\s]$, this implies $\dnfty(M[0,0,\lambda_\s], M[0,1,\lambda_\s], a_1, M[1, 1, \lambda_\s])$.  Since Symmetry for $\s$ holds, there are $M[3, 0, \lambda_\s], M[4, 2, \lambda_\s] \in K_{\lambda_\s}$  such that $M[0,0,\lambda_\s] \prec M[3, 0, \lambda_\s] \prec M[4, 2, \lambda_\s]$ and $M[1, 1, \lambda_\s] \prec M[4, 2, \lambda_\s]$ with $a_1 \in M[3, 0, \lambda_\s]$ and $\dnfty(M[0,0,\lambda_\s], M[3, 0, \lambda_\s], a_2, M[4, 2, \lambda_\s])$.\\
By chasing diagrams, $tp(a_1/ M[0,0,\lambda_\s], M[1, 0, \mu]) = tp(a_1/M[0,0,\lambda_\s], M[3, 0, \lambda_\s])$, so there are $M[3, 0, \mu] \in K_\mu$ and $f:M[1, 0, \mu] \to_{M[0,0,\lambda_\s]} M[3, 0, \mu]$ such that $M[3, 0, \lambda_\s] \prec M[3, 0, \mu]$ and $f(a_1)=a_1$.  
Since $\geq \s$ satisfies Extension Existence and $K$ has the amalgamation property, there is a nonforking extension of $f(tp(a'/M[1, 0, \mu], M[2, 1, \mu]))$ to $M[3, 0, \mu]$.  This means that there are $M[3, 1, \mu] \in K_\mu$, $a'' \in M[3, 1, \mu]$, and $g:M[2, 1, \mu] \to M[3, 1, \mu]$ such that
\begin{itemize}
	\item $M[3, 0, \mu] \prec M[3, 1, \mu]$; 
	\item $f \subset g$;
	\item $tp(a''/f(M[1, 0, \mu]), M[3, 1, \mu]) = tp(g(a')/f(M[1, 0, \mu]), M[3, 1, \mu])$; and
	\item $\dnfty(f(M[1, 0, \mu]), M[3, 0, \mu], a'', M[3, 1, \mu])$
\end{itemize}
Extend $g$ to an $L(K)$-isomorphism $G$ with range including $M[3, 1, \mu]$.

Then $\dnfty(M[1, 0, \mu], G^{-1}(M[3, 0, \mu]), G^{-1}(a''), G^{-1}(M[3, 1, \mu]))$ and $tp(G^{-1}(a'')/M[1, 0, \mu], G^{-1}(M[3, 1, \mu])) = tp(a'/M[1, 0, \mu], M[2, 1, \mu])$.  Since $\dnfty(M[0,0,\lambda_\s], M[1, 0, \mu], a', M[2, 1, \mu])$, this type equality means that $\dnfty(M[0,0,\lambda_\s], M[1, 0, \mu], G^{-1}(a''), G^{-1}(M[3, 1, \mu]))$.  Since $\geq \s$ satisfies Transitivity, we have $\dnfty(M[0,0,\lambda_\s], G^{-1}(M[3, 0, \mu]), G^{-1}(a''), G^{-1}(M[3, 1, \mu]))$.  Since $G \supset g \supset f$ fixes $M[0,0,\lambda_\s]$ and $\geq \s$ satisfies Invariance, we have $\dnfty(M[0,0,\lambda_\s], M[3, 0, \mu], a'', M[3, 1, \mu])$.  By Monotonicity, we have $\dnfty(M[0,0,\lambda_\s], M[3, 0, \lambda_\s], a'', M[3, 1, \mu])$.  Recall that we picked $M[3, 0, \lambda_\s]$ such that \\$\dnfty(M[0,0,\lambda_\s], M[3, 0, \lambda_\s], a_2, M[4, 2, \lambda_\s])$ and that 
\begin{eqnarray*}
tp(a_2/M[0,0,\lambda_\s], M[4, 2, \lambda_\s]) &=& tp(a'/M[0,0,\lambda_\s], M[2, 1, \mu])\\ 
&=& tp(g(a')/M[0,0,\lambda_\s], M[3, 1, \mu])\\
 &=& tp(a''/M[0,0,\lambda_\s], M[3, 1, \mu])
\end{eqnarray*}
since $g$ fixes $M[0,0,\lambda_\s]$.  By Uniqueness, $tp(a_2/ M[3, 0, \lambda_\s], M[4, 2, \lambda_\s]) = tp(a''/ M[3, 0, \lambda_\s], M[3, 1, \mu])$.  Since $a_1 \in M[3, 0, \lambda_\s]$ and $M^- \prec M[0,0,\lambda_\s] \prec M[3, 0, \lambda_\s]$, this implies $tp(a_1 a_2/M^-, M[4, 2, \lambda_\s]) = tp(a_1 a''/M^-, M[3, 1, \mu])$.\\
On the other hand, since $f(a_1) = a_1$ and $f$ fixes $M[0,0,\lambda_\s]$, we have that
\begin{eqnarray*}
tp(a''/f(M[1, 0, \mu]), M[3, 1, \mu]) &=& tp(g(a')/f(M[1, 0, \mu]), M[3, 1, \mu])\\
tp(a_1 a''/f(M[1, 0, \mu]), M[3, 1, \mu]) &=& tp(a_1 g(a')/f(M[1, 0, \mu]), M[3, 1, \mu])\\
tp(a_1 a''/M^-, M[3, 1, \mu]) &=& tp(a_1 g(a')/M^-, M[3, 1, \mu]) = tp(a_1 a'/M^-, M[2, 1, \mu])\\
\end{eqnarray*}
So $tp(a_1 a_2/M^-, M[1, 1, \mu]) = tp(a_1 a'/M^-, M[2, 1, \mu])$, as desired.\\
Since $M^- \prec M[0,0, \mu]$ of size $\lambda_\s$ was arbitrary and $K$ is $\lambda_\s$-tame for 2-types, we have $tp(a_1 a_2/M[0,0, \mu], M[1, 1, \mu]) = tp(a_1 a'/M[0,0, \mu], M[2, 1, \mu])$.  This proves the claim and the theorem.\hfill \dag\\

Thus, we add the following hypothesis.  Note that basic types are only defined for types of length one, so a hypothesis of ``tameness for basic 2-types'' would not make sense.

\begin{hyp} \label{tame2hyp}
$K$ is $\lambda_\s$-tame for 2-types
\end{hyp}

We focus on this method for obtaining Symmetry due to its similarity to Hypothesis \ref{tamehyp}.  However, there is another way to derive Symmetry that does not rely on the structure of extending the frame $\s$.  Recall from Shelah \cite{sh576} that a type $p \in S(M)$ is \emph{minimal} iff it has at most one non-algebraic extension to any $N \succ M$ with $\|N \| = \|M\|$ and that basic types in the frame from Theorem \ref{genframeext} are exactly the rooted minimal types.  Then \cite{shelahaecbook}.II.\S.3.7 combines the minimality of basic types with disjoint amalgamation in $\lambda_\s$ to derive Symmetry for $\s$.  This proof can be adapted to get the following.

\begin{ntheorem}[Without Hypothesis \ref{tame2hyp}]
If basic types for $\s$ are minimal and $K_{\geq \lambda_\s}$ satisfies disjoint amalgamation, then $\geq \s$ satisfies Axiom (E)(f).
\end{ntheorem}

%%%%%%%%%%%%%%%%%%%%%%%%%%%%%%%%%%%%%%%
\section{No Maximal Models} \label{nmmsection}
%%%%%%%%%%%%%%%%%%%%%%%%%%%%%%%%%%%%%%%

Recall that we are working under Hypotheses \ref{basichyp}, \ref{tamehyp}, and \ref{tame2hyp}; these are that $\s$ is a good $\lambda$-frame and $K_{\geq \lambda}$ has amalgamation; that $K$ is $\lambda$-tame for basic 1-types (in the sense of $\geq \s$); and that $K$ is $\lambda$-tame for 2-types.  The results so far have shown that $\geq \s$ is a good frame except possibly for the ``no maximal models'' clause.

In this section, we adapt the proof of \cite{shelahaecbook}.II.\S4.13.3 to show that if $K_{\geq \kappa}$ has a good frame $\geq \s$, then $K_{\geq \kappa}$ has no maximal model.  This is no real change in the proof, except to include the case of where the size of the model is a limit cardinal.  This proof makes use of a strengthening of Non-Forking Amalgamation that Shelah calls Long Non-Forking Amalgamation.  We include a proof of the final result, which combines the work of \cite{shelahaecbook}.II.\S4.9.1, .12.1, and .13.3, to show all of the details.

\begin{ntheorem}[\cite{shelahaecbook}.II.\S4.13.3*] \label{weakstepup}
Assume $\lambda < \kappa$ and $K_\kappa$ is non-empty.  Then $K_\kappa$ has no maximal models.
\end{ntheorem}
{\bf Proof:}  Let $N^0 \in K_\kappa$ and let $\seq{N^0_i \in K_{[\lambda, \kappa)} : i < \kappa}$ be a resolution.  From Density, we know that, for each $i < \kappa$, there is some $a_i \in N^0_{i+1} - N^0_i$ such that $tp(a_i/N^0_{i}, N^0_{i+1}) \in \gS^{bs}_{\geq \s}(N^0_i)$ and some $p \in tp(b/N^0_0, N^1_0) \in \gS^{bs}_{\geq \s}(N^0_0)$; we might have $a_0 = b$ and $N^0_1 = N^1_0$, but this is okay.\\
We will construct, by induction on $\alpha \leq \lambda$, a coherent sequence $\seq{N_\alpha^1, f_{\beta, \alpha}: N_\beta \to N_\alpha \mid \beta < \alpha \leq \lambda}$ such that 
\begin{enumerate}

	\item $N_\alpha^0 \prec N_\alpha^1$;
	\item $\dnf_{\geq \s}(N_\alpha^0, N_{\alpha+1}^0, f_{0, \alpha+1}(b), N_{\alpha+1}^1)$; and 
	\item  $f_0 = \id_{N^1_0}$.

\end{enumerate}
$\alpha = 0$ is already defined.  For $\alpha$ limit, we take a direct limit.  For $\alpha = \beta+1$, we have that $N_\alpha^0 \prec N_\alpha^1, N_{\alpha+1}^0$ with $tp(a_\alpha/N^0_\alpha, N^0_{\alpha+1}), tp(f_{0, \beta}(b)/N^0_\alpha, N^1_\alpha) \in \gS^{bs}_{\geq \s^-}(N^0_\alpha)$.  Then we use Non-Forking Amalgamation to find $f_\beta: N^1_\beta \to N^1_\alpha$ with $N^0_\alpha \prec N^1_\alpha$ so $\dnf_{\geq \s}(N^0_\beta, N^0_{\alpha}, f_\beta(f_{0, \beta}(b)), N^1_{\alpha})$ and $\dnf_{\geq \s}(N^0_\alpha, f_\beta(N^1_\beta), a_\alpha, N^1_\alpha)$.  For $\gamma \leq \beta$, set $f_{\gamma, \alpha} = f_\beta \circ f_{\gamma, \beta} $.\\
This completes our construction.  Now we have that $N^0 = \bigcup_{\alpha < \lambda} N^0_\alpha \precneqq \bigcup_{\alpha < \lambda} N^1_\alpha = N^1 \in K_\lambda$, since $f_{0, \lambda}(b) \not \in N^0$.  Since $N^0 \in K_\lambda$ was arbitrary, we are done. \hfill \dag\\

This allows us to prove the existence of arbitrarily large models.

\begin{ncor} \label{nmm}
$K$ has no maximal models.  In particular, it has models of all cardinalities.
\end{ncor}

%%%%%%%%%%%%%%%%%%%%%%%%%%%%%%%%%%%%%%%
\section{Good Frames} \label{gfsection}
%%%%%%%%%%%%%%%%%%%%%%%%%%%%%%%%%%%%%%%

We drop the previous hypotheses for this section, although $K$ will always be an AEC.

We combine our previous results into the following theorem.
\begin{ntheorem}\label{frameup}
Suppose $K$ is an AEC with amalgamation.  If $K$ has a good $\lambda$-frame $\s$ and is $\lambda_\s$-tame for 1- and 2- types, then $\geq \s$ is a good frame.  
\end{ntheorem}

{\bf Proof:}  From Theorem \ref{basicextension}, we know that $\geq \s$ satisfies all of the axioms of a good frame except for amalgamation, joint embedding, no maximal models, uniqueness, basic stability, extension existence, and symmetry.  Amalgamation and joint embedding follow from the assumption of this theorem.  Uniqueness, basic stability, and extension existence follow from tameness for 1-types by Theorem \ref{uniqtrans}, Corollary \ref{uniqtrans}, and Theorem \ref{exttrans}.  Symmetry follows from tameness for 2-types by Theorem \ref{symtheorem}.  Finally, no maximal models follows from tameness for 1- and 2-types by Corollary \ref{nmm}. \hfill \dag\\

This is the main theorem promised in the introduction.  We provide proofs of some of the other claims as well.  First, we can trade the assumption of no maximal models in the categoricity transfer of \cite{tamenessthree} for a set-theoretic assumption, a slight increase in tameness, and an extra categoricity cardinal.

\begin{ntheorem}\label{above}
Let $K$ be an AEC with amalgamation and $LS(K) < \kappa \leq \lambda$ such that
\begin{enumerate}
	\item $K$ is $\kappa$ tame for 1- and 2-types; and
	\item $K$ is categorical in $\lambda$ and $\lambda^+$ with
	\begin{enumerate}
	\item[$(*)$] $2^\lambda < 2^{\lambda^+} < 2^{\lambda^{++}}$ and $WDmId(\lambda^+)$ is not $\lambda^{++}$-saturated.
\end{enumerate}
\end{enumerate}
Then $K$ is categorical in all $\mu \geq \lambda$.
\end{ntheorem}

{\bf Proof:} By 2. of the hypothesis and Theorem \ref{genframeext}, $K$ has a good $\lambda^+$-frame $\s$.  By Theorem \ref{frameup} and tameness, $\geq \s$ is a good frame.  In particular, $K$ has no maximal models.  Then, we can apply the categoricity transfer of \cite{tamenessthree} to show that $K$ is categorical in all $\mu \geq \lambda^+$ and we have $\mu = \lambda$ as part of the hypothesis. \hfill \dag\\

All in all, this is not a very good trade.  On the other hand, during this proof we constructed our promised independence relation in a tame and categorical AEC.  There are two related sets of assumptions that allow us to do so, both of which utilize the work of Shelah, Grossberg and VanDieren, and Theorem \ref{frameup}.

\begin{nprop}\label{forkingincattame}
Let $K$ be an AEC with amalgamation that is $\kappa$-tame for 1- and 2- types and is categorical in $\lambda^+$ with $\lambda > \kappa > LS(K)$.  If either of the two following hold
\begin{enumerate}

\item $K$ has no maximal models and joint embedding and there is some $\mu \geq \min \{ \lambda^+, \beth_{\chi}\}$ for $\chi = (2^{\beth_{(2^{LS(K)})^+}})^+$ such that $2^\mu < 2^{\mu^+} < 2^{\mu^{++}}$ and $WDmId(\mu^+)$ is not $\mu^{++}$-saturated; or

\item $K$ is categorical in $\lambda$ and $2^\lambda < 2^{\lambda^+} < 2^{\lambda^{++}}$ and $WDmId(\lambda^+)$ is not $\lambda^{++}$-saturated;

\end{enumerate}
then there is a good frame $\geq \s$ with $\lambda_\s = \mu^+$ in case (a) and $\lambda_\s = \lambda^+$ in case (b).
\end{nprop}

{\bf Proof:} Case (b) was handled in Theorem \ref{above} above.  In case (a), the assumption of joint embedding and no maximal models means that we can use the results of \cite{tamenessthree} and \cite{sh394} to conclude that $K$ is categorical in every cardinal above $\min \{ \lambda^+, \beth_{(2^{\beth_{(2^{LS(K)})^+}})^+}\}$; in particular, $\mu$ and $\mu^+$.  Then we can use Theorem \ref{genframeext} to derive a good $\mu^+$ frame $\s$.  By Theorem \ref{frameup}, $\geq \s$ is a good frame with $\lambda_\s = \mu^+$.\hfill \dag\\

%%%%%%%%%%%%%%%%%%%%%%%%%%%%%%%%%%%%%%%
\section{Uniqueness of Limit Models} \label{ulmsection}
%%%%%%%%%%%%%%%%%%%%%%%%%%%%%%%%%%%%%%%

Recall that $M_\alpha$ is a $(\lambda, \alpha)$-limit model over $M_0$ iff there is a continuous, increasing chain $\seq{M_i \in K_\lambda : i \leq \alpha}$ so $M_{i+1}$ is universal over $M_i$ for all $i < \alpha$.  An easy back-and-forth argument shows that a $(\lambda, \theta_1)$-limit model and $(\lambda, \theta_2)$-limit model over $M$ are isomorphic over $M$ if $\cf \theta_1 = \cf \theta_2$.  The general question of uniqueness of limit models asks if this is true for all $\theta_1, \theta_2 < \lambda^+$.  This question is suspected to be very important in the classification theory of AECs and is addressed in Shelah and Villaveces \cite{shvi635}; VanDieren \cite{vandierennomax} \cite{nomaxerrata}; and Grossberg, VanDieren, and Villaveces \cite{gvv}.  An important caveat is that the uniqueness of limit models result of VanDieren \cite{vandierennomax} \cite{nomaxerrata} was born out of a gap she discovered in \cite{shvi635} and works in the context of amalgamation only over unions of limit models, rather than the full amalgamation used here and in \cite{gvv}.  Shelah outlines the proof of the uniqueness of limit models from the existence of a good $\lambda$-frame, culminating in \cite{shelahaecbook}.II.\S4.8.  We fill in the details because the outlines Shelah offers are very sparse (see, for instance, \cite{shelahaecbook}.II.\S4.11) and to hopefully quell the doubts expressed in \cite{gvv}.6.  Primarily, we provide a detailed proof of a weakening of \cite{shelahaecbook}.II.\S4.11 that constructs a matrix of models, the corner of which is both a $(\lambda, \theta_1)$ and $(\lambda, \theta_2)$ limit model over the same base.  

\begin{nlemma}[II.\S4.11-]\label{ulmlemma}
Suppose we have a good $\lambda$-frame $\s$ and
\begin{enumerate}

	\item regular $\theta_1, \theta_2 \leq \lambda$ such that $\delta_1 = \lambda \otimes \theta_1$ and $\delta_2 = \lambda \otimes \theta_2$
	
	\item $M \in K_\lambda$.

\end{enumerate}
Then, we can find functions $\epsilon:\delta_1 \to \delta_2$ and $\eta:\delta_2 \to \delta_1$, an increasing, continuous matrix of models and embeddings $\seq{M_{\alpha, \beta} \in K_\lambda : \alpha \leq \delta_1, \beta \leq \delta_2}$ and coherent $\seq{f_{(\alpha_0, \beta_0)}^{(\alpha_1, \beta_1)}: M_{(\alpha_0, \beta_0)} \to M_{(\alpha_1, \beta_1)} \mid \alpha_0 \leq \alpha_1 \leq \delta_1; \beta_0 \leq \beta_1 \leq \delta_2  }$, and $\seq{b_\alpha^1 \in M_{\alpha+1, \epsilon(\alpha)+1}: \alpha < \delta_1}$ and $\seq{b_\beta^2 \in M_{\eta(\beta)+1, \beta+1} : \beta < \delta_2}$ so

\begin{enumerate}

	\item[$(\gamma)_1$] $tp(f_{(\alpha +1, \epsilon(\alpha)+1)}^{(\alpha+1, \delta_2)}(b_\alpha^1)/f_{(\alpha, \delta_2)}^{(\alpha+1, \delta_2)}(M_{\alpha, \delta_2}), M_{\alpha+1, \delta_2})$ does not fork over $f_{(\alpha, \epsilon(\alpha)+1)}^{(\alpha+1, \delta_2)}(M_{\alpha, \epsilon(\alpha)+1})$.

	\item[$(\gamma)_2$] $tp(f_{(\eta(\beta) +1, \beta+1)}^{(\delta_1, \beta+1)}(b_\beta^2)/f_{(\delta_1, \beta)}^{(\delta_1, \beta+1)}(M_{\delta_1, \beta}), M_{\delta_1, \beta+1})$ does not fork over $f_{(\eta(\beta)+1, \beta)}^{(\delta_1, \beta+1)}(M_{\eta(\beta)+1, \beta})$.

	\item[$(\delta)_1$] For all $\alpha < \delta_1, \beta < \delta_2, p \in \gS^{bs}(M_{\alpha, \beta+1})$, there are $\lambda$ many $\alpha' > \alpha$ such that $\beta = \epsilon(\alpha')$ and $tp(b_{\alpha'}^1 / f^{(\alpha'+1, \beta+1)}_{(\alpha', \beta+1)}(M_{\alpha', \beta+1}), M_{\alpha'+1, \beta+1})$ is a nonforking extension of $f^{(\alpha'+1, \beta+1)}_{(\alpha, \beta+1)}(p)$.

	\item[$(\delta)_2$] For all $\alpha < \delta_1, \beta < \delta_2, p \in \gS^{bs}(M_{\alpha+1, \beta})$, there are $\lambda$ many $\beta' > \alpha$ such that $\alpha = \eta(\beta')$ and $tp(b_{\beta'}^2 / f^{(\alpha+1, \beta'+1)}_{(\alpha+1, \beta')}(M_{\alpha+1, \beta'}), M_{\alpha+1, \beta'+1})$ is a nonforking extension of $f^{(\alpha+1, \beta'+1)}_{(\alpha+1, \beta)}(p)$.

\end{enumerate}
\end{nlemma}

The minus indicates that the original lemma has several clauses that aren't needed for this application, so we drop them.  Our numbering is, again, to be consistent with \cite{shelahaecbook}.  Here, coherent means that for $\alpha_0 \leq \alpha_1 \leq \alpha_2 \leq \delta_1$ and $\beta_0 \leq \beta_1 \leq \beta_2 \leq \delta_2$, we have $f_{(\alpha_0, \beta_0)}^{(\alpha_2, \beta_2)} = f_{(\alpha_1, \beta_1)}^{(\alpha_2, \beta_2)} \circ f_{(\alpha_0, \beta_0)}^{(\alpha_1, \beta_1)}$

{\bf Proof:} There are disjoint $\seq{u^1_{\alpha,1} \subset \delta_1 : \alpha < \delta_1, i < \lambda}$ and $\seq{u^2_{\beta, i} \subset \delta_2: \beta < \delta_2, i < \lambda}$ such that, for each $\ell =1,2$ and each $\alpha, \gamma < \delta_\ell$ and $i < \lambda$, we have
\begin{itemize}
	\item $|u^\ell_{\alpha, i}| = \lambda$; and
	\item $\gamma \in u^\ell_{\alpha, i}$ implies $\gamma > \alpha$.
\end{itemize}
We want to reindex these sequences based on the types of our matrix models to, for instance, $\seq{u^1_{\alpha, \beta, p} \subset \delta_1 : \alpha < \delta_1, \beta < \delta_2, p \in \gS^{bs}(M_{\alpha, \beta+1})}$ by changing the $i$'s to $\beta, p$'s.  Since $|\delta_2| = \lambda$ and $K$ is $bs$-stable in $\lambda$, there is no problem with the cardinalities.  However, we have not defined the models $M_{\alpha, \beta}$ yet.  Formally, we should index these in terms of $\alpha, \beta, j$ for $j < \lambda$ and, once $M_{\alpha, \beta+1}$ is defined, enumerate the types.  However, this adds more complexity to an already technical proof.  Thus, we write them now as $\seq{u^1_{\alpha, \beta, p} \subset \delta_1 : \alpha < \delta_1, \beta < \delta_2, p \in \gS^{bs}(M_{\alpha, \beta+1})}$ and $\seq{u^2_{\alpha, \beta, p} \subset \delta_2 : \alpha < \delta_1, \beta < \delta_2, p \in \gS^{bs}(M_{\alpha+1, \beta})}$, noting that they still satisfy the above properties.
Define $\epsilon: \delta_1 \to \delta_2$ by $\epsilon(\alpha) = \beta$ iff $\alpha \in u^1_{\alpha_0, \beta, p_0}$ and define $\eta:\delta_2 \to \delta_1$ by $\eta(\beta) = \alpha$ iff $\beta \in u^2_{\alpha, \beta_0, p_0}$.  Note that $\epsilon(\alpha) = \beta$ implies $\alpha > \alpha_0$ and $\eta(\beta) = \alpha$ implies $\beta > \beta_0$.

Now we build the rest of our objects by induction so

\begin{enumerate}

	\item $M_{\alpha, 0} = M_{0, \beta} = M$ for all $\alpha \leq \delta_1$ and $\beta \leq \delta_2$.
	
	\item for each $(\alpha, \beta) \in \delta_1 \times \delta_2$, 
	\begin{enumerate}
	
		\item \begin{enumerate}
		
			\item[(i)] if $\epsilon(\alpha) < \beta$, $tp(f^{(\alpha+1, \beta)}_{(\alpha+1, \epsilon(\alpha)+1)}(b_\alpha^1)/f^{(\alpha+1, \beta)}_{(\alpha, \beta)}(M_{\alpha, \beta}), M_{\alpha+1, \beta})$ does not fork over \\$f^{(\alpha+1, \beta)}_{(\alpha, \epsilon(\alpha)+1)}(M_{\alpha, \epsilon(\alpha) + 1})$
			
			\item[(ii)] if $\epsilon(\alpha) = \beta$, then $\alpha \in u^1_{\alpha_0, \beta, p_0}$ for some $\alpha_0 < \alpha$ and $p_0 \in \gS^{bs}(M_{\alpha_0, \beta+1})$ and we pick $b_\alpha^1 \in M_{\alpha+1, \beta+1}$ that realizes the nonforking extension of $f^{(\alpha+1, \beta+1)}_{(\alpha_0, \beta+1)}(p_0)$ to $f^{(\alpha+1, \beta+1)}_{(\alpha, \beta+1)}(M_{\alpha, \beta+1})$.
		
		\end{enumerate}
		
		\item \begin{enumerate}
		
			\item[(i)] if $\eta(\beta) < \alpha$, $tp(f^{(\alpha, \beta+1)}_{(\eta(\beta)+1, \beta+1)}(b_\beta^2)/f^{(\alpha, \beta+1)}_{(\alpha, \beta)}(M_{\alpha, \beta}), M_{\alpha, \beta+1})$ does not fork over \\$f^{(\alpha, \beta+1)}_{(\eta(\beta)+1, \beta)}(M_{\eta(\beta)+1, \beta})$
			
			\item[(ii)] if $\eta(\beta)=\alpha$, then $\beta \in u^2_{\alpha, \beta_0, p_0}$ for some $\beta_0 < \beta$ and $p_0 \in \gS^{bs}(M_{\alpha+1, \beta_0})$ and we pick $b_\beta^2 \in M_{\alpha+1, \beta+1}$ that realizes the nonforking extension of $f^{(\alpha+1, \beta+1)}_{(\alpha+1, \beta_0)}(p_0)$ to $f^{(\alpha+1, \beta+1)}_{(\alpha+1, \beta)}(M_{\alpha+1, \beta})$.
		
		\end{enumerate}
	
	\end{enumerate}

\end{enumerate}

{\bf Construction:} The edges of the matrices are our base cases.\\
If $\alpha$ or $\beta$ is limit, then we construct the model via direct unions and check that our conditions hold.\\
So we are in the case where we have $\alpha < \delta_1$ and $\beta < \delta_2$ and we need to construct $M_{\alpha + 1, \beta+1}$ and the embeddings given $M_{\alpha +1, \beta}$ and $M_{\alpha, \beta+1}$.  Before we construct our model, we do some preparatory work and find $N_\alpha \succ M_{\alpha, \beta+1}$ and $N_\beta \succ M_{\alpha+1, \beta}$; $a_\alpha \in N_\alpha - M_{\alpha, \beta+1}$ and $a_\beta \in N_\beta -  M_{\alpha+1, \beta}$; and $n_\alpha \in M_{\alpha+1, \beta}$ so its type over $f_{(\alpha, \beta)}^{(\alpha+1, \beta)}(M_{\alpha, \beta})$ is basic and $n_\beta \in M_{\alpha, \beta+1}$ so its type over $f_{(\alpha, \beta)}^{(\alpha, \beta+1)}(M_{\alpha, \beta})$ is basic.\\
\begin{enumerate}

	\item If $\epsilon(\alpha) < \beta$, then we have $tp(f^{(\alpha+1, \beta)}_{(\alpha+1, \epsilon(\alpha)+1)}(b_\alpha^1) / f^{(\alpha+1, \beta)}_{(\alpha, \beta)}(M_{\alpha, \beta}, M_{\alpha+1, \beta})$ is basic, so pick $n_\alpha = f^{(\alpha+1, \beta)}_{(\alpha+1, \epsilon(\alpha)+1)}(b^1_\alpha)$.  Otherwise, use the Density to pick $n_\alpha$ arbitrarily.  Note that this axiom is not necessary, but helps to make our construction more symmetric.
	
	\item If $\epsilon(\alpha) = \beta$, then $\alpha \in u^1_{\alpha_0, \beta, p_0}$ by construction, so by Extension Existence, there is $tp(a_\alpha/ M_{\alpha, \beta+1}, N_\alpha)$ that is a nonforking extension of $f_{(\alpha_0, \beta+1)}^{(\alpha, \beta+1)}(p_0)$.  Otherwise, pick them arbitrarily.  Note that $\alpha_0 < \alpha$, so $M_{\alpha_0, \beta+1}$ has been constructed prior to this step, so this enumeration is well defined.
	
	\item If $\eta(\beta) < \alpha$, then we have $tp(f^{(\alpha, \beta+1)}_{(\eta(\beta)+1, \beta+1)}(b_\beta^2)/f_{(\alpha, \beta)}^{(\alpha, \beta+1)}(M_{\alpha, \beta}), M_{\alpha, \beta+1})$ is basic, so pick $n_\beta = f^{(\alpha, \beta+1)}_{(\eta(\beta)+1, \beta+1)}(b^2_\beta)$.  Otherwise, pick $n_\beta$ arbitrarily.
	
	\item If $\eta(\beta) = \alpha$, then $\beta \in u^2_{\alpha, \beta_0, p_0}$, so find, by Extension Existence, $tp(a_\beta/M_{\alpha+1, \beta}, N_\beta)$ that is a nonforking extension of $f_{(\alpha+1, \beta_0)}^{(\alpha+1, \beta)}(p_0)$.  Otherwise, pick them arbitrarily.  As above, $\beta_0 < \beta$, so this is well defined.

\end{enumerate}

Now that we have this, we apply Non-Forking Amalgamation to get the following
\[
\xymatrix{
N_\beta \ar[rr]^{g_\beta} & & M_{\alpha+1, \beta+1}\\
M_{\alpha+1, \beta} \ar[u] & & \\
M_{\alpha, \beta} \ar[u]_{f_{(\alpha, \beta)}^{(\alpha+1, \beta)}} \ar[r]_{f_{(\alpha, \beta)}^{(\alpha, \beta+1)}} & M_{\alpha, \beta+1} \ar[r] & N_\alpha \ar[uu]_{g_\alpha}
}
\]
Set $f^{(\alpha+1, \beta+1)}_{(\alpha+1, \beta)} = g_\beta \rest M_{\alpha+1, \beta}$ and $f_{(\alpha, \beta+1)}^{(\alpha+1, \beta+1)} = g_\alpha \rest M_{\alpha, \beta+1}$.  Then compose the rest of the embeddings to make everything coherent.
\begin{enumerate}

	\item If $\epsilon(\alpha) < \beta$, then $\epsilon(\alpha) < \beta+1$ and nonforking amalgamation tells us (after a little rewriting) that 
	\begin{eqnarray}
	tp(f_{(\alpha + 1, \epsilon(\alpha)+1)}^{(\alpha+1, \beta+1)}(b_\alpha+1) / f_{(\alpha, \beta+1)}^{(\alpha+1, \beta+1)}(M_{\alpha, \beta+1}), M_{\alpha+1, \beta+1})\notag\\
	\te{does not fork over }f_{(\alpha, \beta)}^{(\alpha+1, \beta+1)}(M_{\alpha, \beta})
	\end{eqnarray}  
	By induction, we have that $tp(f^{(\alpha+1, \beta)}_{(\alpha+1, \epsilon(\alpha)+1)}(b_\alpha^1)/f^{(\alpha+1, \beta)}_{(\alpha, \beta)}(M_{\alpha, \beta}), M_{\alpha+1, \beta})$ does not fork over $f^{(\alpha+1, \beta)}_{(\alpha, \epsilon(\alpha)+1)}(M_{\alpha, \epsilon(\alpha) + 1})$.  Applying $f_{(\alpha+1, \beta)}^{(\alpha+1, \beta+1)}$ to this and applying Monotonicity, we get that 
	\begin{eqnarray}
	tp(f^{(\alpha+1, \beta+1)}_{(\alpha+1, \epsilon(\alpha)+1)}(b_\alpha^1)/f^{(\alpha+1, \beta+1)}_{(\alpha, \beta)}(M_{\alpha, \beta}), M_{\alpha+1, \beta+1})\notag\\
	\te{ does not fork over }f^{(\alpha+1, \beta+1)}_{(\alpha, \epsilon(\alpha)+1)}(M_{\alpha, \epsilon(\alpha) + 1})
	\end{eqnarray}
	Then, we apply Transitivity to Eqs. (9.1) and (9.2) and get that \\$tp(f_{(\alpha + 1, \epsilon(\alpha)+1)}^{(\alpha+1, \beta+1)}(b_{\alpha+1}) /$ $ f_{(\alpha, \beta+1)}^{(\alpha+1, \beta+1)}(M_{\alpha, \beta+1}), M_{\alpha+1, \beta+1})$ does not fork over $f^{(\alpha+1, \beta+1)}_{(\alpha, \epsilon(\alpha)+1)}$ $(M_{\alpha, \epsilon(\alpha) + 1})$, as desired.
 
	\item If $\epsilon(\alpha) = \beta$, then we set $b^1_\alpha = g_\alpha(a_\alpha) \in M_{\alpha+1, \beta+1}$.  We know that $a_\alpha \realize f_{(\alpha_0, \beta+1)}^{(\alpha, \beta+1)}(p_0)$, so $b^1_\alpha$ realizes $f_{(\alpha_0, \beta+1)}^{(\alpha+1, \beta+1)}(p_0)$.  Additionally, we picked $a_\alpha$ so $tp(a_\alpha/ M_{\alpha, \beta+1}, N_\alpha)$ does not fork over $f_{(\alpha_0, \beta+1)}^{(\alpha, \beta+1)}(M_{\alpha, \beta+1})$.  Applying $g_\alpha \supset f_{(\alpha, \beta+1)}^{(\alpha+1, \beta+1)}$ to this and using Monotonicity, we get that $tp(b^1_\alpha/ f_{(\alpha, \beta+1)}^{(\alpha+1, \beta+1)}(M_{\alpha, \beta+1}),$ $ M_{\alpha +1, \beta+1})$ does not fork over $f_{(\alpha_0, \beta+1)}^{(\alpha+1, \beta+1)}(M_{\alpha_0, \beta+1})$.

	\item If $\eta(\beta) < \alpha$ or $\eta(\beta) = \alpha$, the proof is symmetric, since our goal and our set-up is symmetric.
	
\end{enumerate}
{\bf This is enough:} Now we want to show that our construction has fulfilled the lemma.
\begin{enumerate}

	\item[$(\gamma)_1$]  {Set $\alpha < \delta_1$.}  For each $\beta > \epsilon(\alpha)$, we know that $tp(f^{(\alpha+1, \beta)}_{(\alpha+1, \epsilon(\alpha)+1)}(b_\alpha^1)/f^{(\alpha+1, \beta)}_{(\alpha, \beta)}(M_{\alpha, \beta}),$ $ M_{\alpha+1, \beta})$ does not fork over $f^{(\alpha+1, \beta)}_{(\alpha, \epsilon(\alpha)+1)}(M_{\alpha, \epsilon(\alpha) + 1})$ by 2.(a)(i) of the construction.  If we apply the map $f_{(\alpha+1, \beta)}^{(\alpha+1, \delta_2)}$ and use Monotonicity, we get that $tp(f^{(\alpha+1, \delta_2)}_{(\alpha+1, \epsilon(\alpha)+1)}(b_\alpha^1)/f^{(\alpha+1, \delta_2)}_{(\alpha, \beta)}(M_{\alpha, \beta}), M_{\alpha+1, \delta_2})$ does not fork over $f^{(\alpha+1, \delta_2)}_{(\alpha, \epsilon(\alpha)+1)}(M_{\alpha, \epsilon(\alpha) + 1})$ for every $\epsilon(\alpha) < \beta < \delta_2$.  Then, by Continuity, we have that $tp(f^{(\alpha+1, \delta_2)}_{(\alpha+1, \epsilon(\alpha)+1)}(b_\alpha^1)/$ $f^{(\alpha+1, \delta_2)}_{(\alpha, \delta_2)}(M_{\alpha, \delta_2}), M_{\alpha+1, \delta_2})$ does not fork over $f^{(\alpha+1, \delta_2)}_{(\alpha, \epsilon(\alpha)+1)}(M_{\alpha, \epsilon(\alpha) + 1})$, as desired.
	
	\item[$(\delta)_1$] Fix $\alpha < \delta_1, \beta < \delta_2, p \in \gS^{bs}(M_{\alpha, \beta+1})$.  Then $u^1_{\alpha, \beta, p} = \epsilon^{-1}(\{\beta\})$ has size $\lambda$ and, for every such $\alpha'$, $tp(b_{\alpha'}^1 / f^{(\alpha'+1, \beta+1)}_{(\alpha', \beta+1)}(M_{\alpha', \beta+1}), M_{\alpha'+1, \beta+1})$ is a nonforking extension of $f^{(\alpha'+1, \beta+1)}_{(\alpha, \beta+1)}(p)$ by 2.(a).(oo).
	
	\item[$(\gamma)_2, (\delta)_2$]  Similarly.

\end{enumerate}
This completes the proof of the lemma. \hfill \dag\\

For reference and, in particular,for use in Boney and Grossberg \cite{shorttamedep}, we note that the only frame properties used were Amalgamation, Density, $bs$-stability, Monotonicity, Transitivity, Symmetry, Extension Existence, and Continuity.  In particular, Continuity was only used for chains of length $\theta_1$ and $\theta_2$.  We can now prove the uniqueness of limit models.

\begin{ntheorem}[\cite{shelahaecbook}II.\S4.8]
If we have a good $\lambda$-frame, then $K_\lambda$ has unique limit models.
\end{ntheorem}

{\bf Proof:}  Let $N_1$ be a $(\lambda, \theta_1)$-limit model over $M$ and $N_2$ be a $(\lambda, \theta_2)$-limit model over $M$.  Apply the lemma above to get functions $\epsilon:\delta_1 \to \delta_2$ and $\eta:\delta_2 \to \delta_1$ and an increasing, continuous matrix of models and embeddings $\seq{M_{\alpha, \beta} \in K_\lambda : \alpha \leq \delta_1, \beta \leq \delta_2}$ and coherent $\seq{f_{(\alpha_0, \beta_0)}^{(\alpha_1, \beta_1)}: M_{(\alpha_0, \beta_0)} \to M_{(\alpha_1, \beta_1)} \mid \alpha_0 \leq \alpha_1 \leq \delta_1; \beta_0 \leq \beta_1 \leq \delta_2  }$ and $\seq{b_\alpha^1 \in M_{\alpha+1, \epsilon(\alpha)+1}: \alpha < \delta_1}$ and $\seq{b_\beta^2 \in M_{\eta(\beta)+1, \beta+1} : \beta < \delta_2}$ as there.\\
By renaming, we get increasing continuous $\seq{M^{\delta_2}_\alpha : \alpha \leq \delta_1}$ and $\seq{M^{\delta_1}_{\beta} : \beta \leq \delta_2}$ such that $M^{\delta_2}_0 = M^{\delta_1}_0 = M$ and $M^{\delta_1}_{\delta_2} = M^{\delta_2}_{\delta_1}$, which is the renaming of $M_{\delta_1, \delta_2}$ with the property
\begin{enumerate}

	\item[$(*)_1$] if $\alpha < \delta_1$ and $p \in \gS^{bs}(M^{\delta_2}_\alpha)$, then there are $\lambda$-many $\alpha' > \alpha$ such that $tp(b^1_{\alpha'}/M^{\delta_2}_{\alpha'}, M^{\delta_2}_{\alpha'+1})$ is a nonforking extension of $p$.
	
	\item[$(*)_2$] if $\beta < \delta_2$ and $p \in \gS^{bs}(M^{\delta_1}_\beta)$, there there are $\lambda$-many $\beta' > \beta$ such that $tp(b^2_{\beta'}/M^{\delta_2}_{\beta'}, M^{\delta_2}_{\beta'+1})$ is a nonforking extension of $p$.

\end{enumerate}

Once we have established these, we use \cite{shelahaecbook}.II.\S4.3 (see Theorem \ref{below}) to see that $M^{\delta_2}_{\delta_1}$ is $(\lambda, \theta_1)$-limit over $M$ and $M^{\delta_1}_{\delta_2}$ is $(\lambda, \theta_2)$-limit over $M$.  Then, by uniqueness of limit models of the same length, we get that 
$$N_1 \cong_M M^{\delta_1}_{\delta_1} = M^{\delta_1}_{\delta_2} \cong_M N_2$$
\hfill \dag\\

For reference, \cite{shelahaecbook}.II.\S.4.3 is stated below and has a detailed proof at the reference and uses only Density and Local Character.

\begin{ntheorem}[Shelah]\label{below}
Assume $\s$ is a good $\lambda$-frame and
\begin{enumerate}

	\item $\delta < \lambda^+$ is a limit ordinal divisible by $\lambda$;
	
	\item $\seq{M_\alpha \in K_\lambda : \alpha \leq \delta}$ is increasing and continuous; and
	
	\item if $i < \delta$ and $p \in \gS^{bs}_\s(M_i)$, then for $\lambda$-many ordinals $j \in (i, \delta)$, there is $c \in M_{j+1}$ realizing the nonforking extension of $p$ in $\gS^{bs}(M_j)$.

\end{enumerate}
Then $M_\delta$ is $(\lambda, \cf \delta)$-limit over $M_0$ and (therefore) universal over it.
\end{ntheorem}

\section{Good Frames in Hart-Shelah} \label{hashframesection} \label{framehashsection}

In this section, we show that some additional hypothesis is necessary to extend a good $\lambda$-frame $\s$ to a good frame $\geq \s$.  This example  was included in response to a referee question about Theorem \ref{basicextension}, and I would like to thank the referee for the question and Alexei Kolesnikov for helpful discussions.

We recall the main result from \cite{untame}.
\begin{ntheorem}[\cite{untame}]
For each $n < \omega$, there is $\phi_n \in L_{\omega_1, \omega}$ so
\begin{enumerate}
	\item $\phi_n$ is categorical in all $\mu \leq \aleph_n$;
	\item $\phi_n$ is not $\aleph_n$-stable;
	\item $\phi_n$ is not categorical in any $\mu > \aleph_n$;
	\item $\phi_n$ has the disjoint amalgamation property; and
	\item if $n > 0$, then
	\begin{enumerate}
		\item $\phi_n$ is $(\aleph_0, \aleph_{n-1})$-tame; in fact, Galois types over models of size $\leq \aleph_{n-1}$ are first order, syntactic types;
		\item $\phi_n$ is $\mu$-stable for $\mu < \aleph_n$; and
		\item $\phi_n$ is not $(\aleph_{n-1}, \aleph_n)$-tame.
	\end{enumerate}
\end{enumerate}
\end{ntheorem}
Note that the sentences $\phi_n$ have been reindexed (as compared to \cite{untame}) in order to avoid unnecessary subscripts such that ``$\phi_n$'' here is ``$\phi_{n+2}$'' there.  We will not give the full definition of $\phi_n$ (it can be found in \cite{untame}.\S1), but will outline some of the key features.  Each model $M$ consists of an index set $I(M)$ (often called the spine) and additional elements built off of this spine, mainly variously indexed copies of $\mathbb{Z}_2$ including fibers over $[I(M)]^{n+2}$ consisting of elements of from the direct sum of $\mathbb{Z}_2$ indexed by $[I(M)]^{n+2}$.  Included in the language are also various projection functions and addition functions.   Added to this is an $(n+3)$-ary predicate $Q$ which codes the addition of $n+2$ many fibers without explicitly including it.

\cite{untame} improves on (and introduces a minor correction to) the original analysis in \cite{hash323}.  In addition to the theorem above, they show that the class of models of $\phi_n$ is model complete (\cite{untame}.4.8).

If $n > 0$, then $\phi_n$ is categorical in at least two successive cardinals ($\aleph_n$ and $\aleph_{n-1}$, for instance), so the results of Shelah \cite{sh576} imply that there is a good $\lambda$-frame under favorable cardinal arithmetic (recall Theorem \ref{genframeext}).  However, the Hart-Shelah example is well-enough understood that cardinal arithmetic is not needed for the existence of a good $\lambda$-frame in this case.  Additionally, we have the existence of a good $\aleph_0$-frame in $\phi_0$, which is only categorical in $\aleph_0$, a result not predicted by \cite{shelahaecbook}.II.

\begin{ntheorem}
Fix $n < \omega$ and $\mu < \aleph_n$.  There is $\s^n_\mu$ such that
\begin{enumerate}
	\item $\s^n_\mu$ is a good $\mu$-frame for $\phi_n$;  
	\item if $\mu < \mu' < \aleph_n$, then $(\geq \s^n_{\mu}) \rest \mu ' = \s^n_{\mu'}$; and
	\item if $\mu' \geq \aleph_n$, then $(\geq \s^n_{\mu}) \rest \mu'$ is a good $\mu$-frame for $\phi_n$ \emph{except} for Uniqueness and Basic Stability, both of which fail.
\end{enumerate}
\end{ntheorem}

Although this proof does not assume any cardinal arithmetic and, therefore, does not use the results of \cite{sh576} to find a frame, the frame definition given is inspired by that frame.

{\bf Proof:}  Fix $n < \omega$.  Then, for this proof, we set $K^n$ to be the models of $\phi_n$ from \cite{untame} and set $M \prec^n N$ iff $M \prec_{L_{\omega_1, \omega}} N$.  This is the same as $M \subset N$ by model completeness.

Fix $\mu < \aleph_n$.  We define the frame $\s^n_\mu = (K^n_\mu, \dnf_\mu^n, \gS^{\bs}_{\mu, n})$ by:
\begin{itemize}
	\item for $M \in K^n$, $tp(a/M, N) \in \gS^{bs}_{\mu, n}(M)$ iff $a \in I(M) - I(N)$; and
	\item $\dnf_\mu^n(M_0, M_1, a, M_3)$ iff $M_0 \prec^n M_1 \prec^n M_3$ and $a \in I(M_3) - I(M_1)$.
\end{itemize}

From the definitions, it then follows that, for any $M \in K_{\geq \mu}$,
\begin{itemize}

	\item $tp(a/M, N) \in \gS^{\bs}_{\geq (\s^n_\mu)}(M)$ iff $a \in I(M) - I(N)$; and
	
	\item $\dnf_{\geq (\s^n_\mu)}(M_0, M_1, a, M_3)$ iff $M_0 \prec^n M_1 \prec^n M_3$ and $a \in I(M_3) - I(M_1)$.

\end{itemize}
This establishes 2.  To show 1. and 3., we will show that $\geq(s^n_\mu)$ satisfies all of the good frame axioms except $\bs$-Stability and Uniqueness and that $\bs$-Stability and Uniqueness hold if the models are of size $< \aleph_n$.  We do this by going through the axioms of Definition \ref{framedef} and showing that they hold.  For notational ease, set $K:=K^n$, $\s:=\s^n_\mu$, $\dnf = \dnf_\mu^n$, and $\gS^{\bs}:= \gS^{\bs}_{\mu, n}$.  Many of the frame properties follow immediately from the definition and the observation that, given $p \in \gS^{\bs}(M)$ and $M_0 \prec M$, $p$ does not fork over $M_0$.  The non-trivial arguments are given below.

\begin{enumerate}
			
	\item[(C)] By \cite{untame}.3.1, $K$ has the stronger property of disjoint amalgamation.  By \cite{untame}.2.15, $K$ is categorical in $\aleph_0$.  Combining this with amalgamation implies that $K$ has joint embedding.  We know that $K$ has arbitrarily large models by \cite{untame}.1.3.  This, plus amalgamation and joint embedding from above, show $K$ has no maximal models; see \cite{bls1003}.3.3.  Thus, $K_\mu$ has no maximal models.  This can also be seen directly be extending the spine, $I$.
	
	\item[(D)] \begin{enumerate}

		\item[(c)] {\bf Density:} The elements of $M$ are determined (up to isomorphism) by $I(M)$.  Thus, $M \precneqq N$ implies $I(M) \subsetneq I(N)$.
		
		\item[(d)] {\bf $bs$-stability:} Below $\aleph_n$, full stability holds by \cite{untame}.7.1; this clearly implies $\bs$-stability.  At $\aleph_n$ and above, the proof of \cite{untame}.6.1 show that there are the maximal number of Galois types of elements from $I$.
	
	\end{enumerate}

	\item[(E)] \begin{enumerate}
	
		\item[(e)] {\bf Uniqueness:} By \cite{untame}.5.1 , Galois types of finite tuples over models of size less than $\aleph_n$ are syntactic, first-order types.  Any two non-algebraic elements in the spine have the same syntactic type, so Uniqueness holds.  At $\aleph_n$ and above, the proof of \cite{untame}.6.8 shows that tameness for basic types fails, so, by Theorem \ref{uniqtrans}, Uniqueness fails as well.
	
		\item[(f)] {\bf Symmetry:} Let $M_0 \prec M_1 \prec M_3$ with $a_1 \in I(M_1) - I(M_0)$ and $a_2 \in I(M_3) - I(M_1)$.  Take $M_2$ to be the substructure generated by $M_0$ and $a_2$ in $M_3$.  Then $I(M_2) = I(M_0) \cup \{ a_2 \}$ and, in particular, $a_1 \not \in I(M_2)$, as desired.
	
		\item[(g)] {\bf Extension Existence:} Let $M$ and $p \in \gS^{\bs}_{\geq \s}(M)$ and $N \succ M$.  Set $p = tp(a/M, N')$ and find a disjoint amalgam $N^* \succ N$ and $f: N' \to_M N^*$.  Then $q = tp(f(a)/N, N^*)$ is a nonforking extension of $p$. \hfill \dag

	\end{enumerate}

\end{enumerate}

In addition to showing that some additional hypothesis is needed to extend a good frame, this example gives a non-trivial example of a frame in ZFC, i.e. without cardinal arithmetic assumptions.  Additionally, this gives an example of a partially categorical AEC with a supersimple-like independence notion, that is, one that has Local Character, Extension Existence, etc., but not Uniqueness.

\end{document}